\documentclass[11pt]{article}
\usepackage{amsfonts}
\usepackage{color}
\newtheorem{theo}{Theorem}[section]
\newtheorem{lem}[theo]{Lemma}
\newtheorem{prop}[theo]{Proposition}

\newcommand{\mysection}[1]{\section{#1} \setcounter{equation}{0}}
\newcommand{\proof}{{\sc Proof.} \quad}
\newcommand{\proofc}{{\sc Proof} \ }
\newcommand{\be}{\begin{equation} \label}
\newcommand{\ee}{\end{equation}}
\newcommand{\bea}{\begin{eqnarray}\label}
\newcommand{\eea}{\end{eqnarray}}
\newcommand{\bas}{\begin{eqnarray*}}
\newcommand{\eas}{\end{eqnarray*}}
\newcommand{\bit}{\begin{itemize}}
\newcommand{\eit}{\end{itemize}}
\newcommand{\qed}{\hfill$\Box$ \vskip.2cm}
\newcommand{\nn}{\nonumber}
\newcommand{\R}{\mathbb{R}}

\newcommand{\pO}{\partial\Omega}

\newcommand{\eps}{\varepsilon}

\newcommand{\supp}{{\rm supp} \, }

\newcommand{\io}{\int_\Omega}
\newcommand{\bom}{\overline{\Omega}}
\newcommand{\mint}{- \hspace*{-4mm} \int}
\newcommand{\abs}{\\[5pt]}
\newcommand{\tm}{T_{max}}
\newcommand{\kf}{k_f}
\newcommand{\Kf}{K_f}
\newcommand{\uw}{\underline{w}}
\begin{document}
\enlargethispage{10mm}
\title{A critical blow-up exponent for flux limitation\\ in a Keller-Segel system}
\author{
Michael Winkler\footnote{michael.winkler@math.uni-paderborn.de}\\
{\small Institut f\"ur Mathematik, Universit\"at Paderborn,}\\
{\small 33098 Paderborn, Germany} }
\date{}
\maketitle
\begin{abstract}
\noindent 
  The parabolic-elliptic cross-diffusion system 
  \bas
    	\left\{ \begin{array}{l}
	u_t = \Delta u - \nabla \cdot \Big(uf(|\nabla v|^2) \nabla v \Big), \\[1mm]
	0 = \Delta v - \mu + u,
	\qquad \io v=0,
	\qquad
	\mu:=\frac{1}{|\Omega|} \io u dx,
 	\end{array} \right.
  \eas
  is considered along with homogeneous Neumann-type boundary conditions in a smoothly bounded domain $\Omega\subset\R^n$, $n\ge 1$,
  where $f$ generalizes the prototype given by
  \bas
	f(\xi) = (1+\xi)^{-\alpha},
	\qquad \xi\ge 0,
	\qquad \mbox{for all } \xi\ge 0,
  \eas
  with $\alpha\in\R$.\abs
  In this framework, the main results assert that if $n\ge 2$, $\Omega$ is a ball and 
  \bas
	\alpha<\frac{n-2}{2(n-1)},
  \eas
  then throughout a considerably large set of radially symmetric initial data, an associated initial value problem
  admits solutions blowing up in finite time with respect to the $L^\infty$ norm of their first components.\abs
  This is complemented by a second statement which ensures that in general and not necessarily symmetric settings,
  if either $n=1$ and $\alpha\in\R$ is arbitrary, or $n\ge 2$ and $\alpha>\frac{n-2}{2(n-1)}$,
  then any explosion is ruled out in the sense that for arbitrary nonnegative and continuous initial data,
  a global bounded classical solution exists.\abs
\end{abstract}
\newpage
\section{Introduction}\label{intro}
We consider nonnegative solutions of the parabolic-elliptic cross-diffusion system
\be{00}
    	\left\{ \begin{array}{l}
	u_t = \Delta u - \nabla \cdot \Big(uf(|\nabla v|^2) \nabla v \Big), \\[1mm]
	0 = \Delta v - \mu + u,
	\qquad \io v=0,
	\qquad 
	\mu:=\frac{1}{|\Omega|} \io u dx,
 	\end{array} \right.
\ee
where the given function $f$ appropriately generalizes the prototype determined by
\be{pr}
	f(\xi) = (1+\xi)^{-\alpha},
	\qquad \xi\ge 0,
	\qquad \mbox{for all } \xi\ge 0,
\ee
with $\alpha\in\R$.
In mathematical biology, systems of this form arise as refinements of the classical Keller-Segel
model for chemotaxis processes,
that is, for processes during which individuals within a population, represented through its density $u=u(x,t)$,
partially orient their diffusive movement toward increasing concentrations $v=v(x,t)$ of a chemical substance produced 
by themselves (\cite{KS}).\abs
Deviating from the simplest choice $f\equiv 1$ that underlies the minimal and most thoroughly studied representative within the
class of Keller-Segel systems (\cite{KS}, \cite{hillen_painter2009}), (\ref{00}) with more general $f$ accounts for refined
approaches in the more recent modeling literature, according to which in several relevant application frameworks
an adequate description of chemotactic motion should include certain gradient-dependent limitations of cross-diffusive fluxes
in the style of (\ref{00})-(\ref{pr}) (\cite{bellomo_flim}, \cite{perthame}, 
\cite{bianchi_painter_sherratt}, \cite{bianchi_painter_sherratt2016}).\abs
In light of well-known results asserting the occurrence of finite-time blow-up in the case when $f\equiv 1$
(\cite{JL}, \cite{nagai1995} , \cite{nagai2001}, \cite{biler}, \cite{suzuki_book}), a natural question seems to be how far
(\ref{00}) retains a certain ability to support singularity formation also in cases when $f$ reflects some saturation effects 
al large signal gradients in the sense that $f(\xi)\to 0$ as $\xi\to\infty$.
Since it can readily be seen that no such explosion arises when $f(\xi)=(1+\xi)^{-\frac{1}{2}}$ for $\xi\ge 0$
(cf.~also Proposition \ref{prop19}), in the context of (\ref{00})-(\ref{pr}) and in a slightly more ambitious formulation,
this amounts to locating the number
\be{ac}
	\alpha_c:=\sup \Big\{ \alpha \in \R \ \Big| \ 
	\mbox{(\ref{00})-(\ref{pr}) admits at least one solution blowing up in finite time} \Big\},
\ee
and hence to deciding, in dependence on whether or not $\alpha_c$ belongs to the open interval $(0,\frac{1}{2})$, if $\alpha_c$
plays the role of a critical exponent that corresponds to a genuinely critical nonlinearity in (\ref{00}).\abs
{\bf The challenge of detecting critical parameter settings in Keller-Segel systems.}\quad
Here we remark that in the context of Keller-Segel type systems, the identification of explosion-critical constellations has 
successfully been accomplished only in quite a small number of cases yet, which may be viewed as reflecting an apparent lack of
appropriate methods for blow-up detection in such classes of cross-diffusion problems.
Indeed, the literature provides a rich variety of techniques capable of discovering situations in which the respective
dissipative ingredients overbalance cross-diffusive destabilization and hence blow-up is ruled out (cf.~\cite{horstmann_DMV} and
\cite{BBTW} for an incomplete overview).
Only in a relatively small number of cases, however, analytical approaches are available which allow for findings on 
singularity formation in parameter settings complementary to those for which results on global existence and boundedness are
available.
Comparatively
well-understood in this respect are thus only certain particular relatives of (\ref{00}) which share some essential 
features with simpler classes of scalar parabolic problems.\abs
Quite far-reaching results are available, e.g., for variants of (\ref{00}) in which the respective chemotactic sensitivity function
depends on the population density $u$ itself, rather than its gradient, such as in the framework of quasilinear systems
with their corresponding first equation given by
\be{Q}
	u_t= \nabla \cdot (D(u)\nabla u) - \nabla \cdot (uS(u)\nabla v).
\ee
Namely, when supplemented either by elliptic equations of the form in (\ref{00}), or even by more complex fully parabolic equations 
for $v$, such modifications of (\ref{00}) admit favorable gradient-like structures that provide accessibility to energy-based
arguments both in the development of global existence theories and in the derivation of blow-up results, and accordingly a
fairly comprehensive knowledge concerning the emergence of singularities could be achieved:
In the prototypical context determined by the choices $D(\xi)=(\xi+1)^{m-1}$ and $S(\xi)=(\xi+1)^{q-1}$, $\xi\ge 0$, for instance,
associated Neumann problems in bounded domains $\Omega\subset\R^n$ admit global bounded solutions for widely arbitrary 
initial data if $m\in\R$ and $q\in\R$ satisfy $m-q>\frac{n-2}{n}$ (\cite{taowin_subcrit}, \cite{djie_win},
\cite{senba_suzuki_AAA}, \cite{kowalczyk_szymanska}),
whereas unbounded solutions exist when $\Omega$ is a ball and $m-q<\frac{n-2}{n}$
(\cite{cieslak_stinner_JDE2012}, \cite{cieslak_stinner_JDE2015}, \cite{cieslak_laurencot_DCDS},
\cite{win_collapse},
\cite{djie_win}; cf.~also \cite{cieslak_laurencot_ANIHPC}).
In some subcases of the latter, the use of Lyapunov functionals even allowed for the construction of certain global solutions
which blow up in infinite time 
(\cite{cieslak_stinner_JDE2012}, \cite{cieslak_stinner_JDE2015}, \cite{lankeit_ITBU}, \cite{win_ITBU}).\abs
Beyond this, however, most classes of relevant chemotaxis systems appear to lack comparable energy structures, and accordingly
the few further studies concerned with rigorous blow-up detection rely on adequately designed ad hoc methods.
In consequence, the identification of explosion-critical parameter settings has so far been achieved only in a small number
of additional cases, in most of which 
either the derivation of collapsing ordinary differential inequalities for moment-like functionals in the flavor
of classical blow-up proofs for semilinear heat equations (\cite{kaplan}, \cite{deng_levine}), or 
even a reduction to single parabolic equations allowing for comparison with exploding subsolutions is possible;
examples of this flavor address critical exponents in nonlinear 
signal production rates (\cite{liu_tao}, \cite{win_NON_ct_signal_critexp}),
optimal conditions on zero-order degradation with respect to blow-up (\cite{fuest_logistic}, \cite{tello_win}),
threshold mass levels for singularity formation (\cite{JL}, \cite{NSY}, \cite{taowin_JEMS}), or also critical interplay
of several ingredients (\cite{bellomo_win_TRAN}).\abs
{\bf Main results.} \quad
The present work now intends to address the corresponding issue of criticality in (\ref{00}) by a {\em combination}
of an essentially moment-based approach with a comparison argument in a crucial first step of an associated blow-up argument.
To substantiate this in the context of a full initial-boundary value problem, let us henceforth consider
\be{0}
    	\left\{ \begin{array}{rcll}
	u_t &=& \Delta u - \nabla \cdot \Big(uf(|\nabla v|^2) \nabla v \Big),
	\qquad & x\in\Omega, \ t>0, \\[1mm]
	0 &=& \Delta v - \mu + u,
	\qquad \io v=0,
	\qquad & x\in\Omega, \ t>0, \\[1mm]
	& & \hspace*{-13mm}
	\frac{\partial u}{\partial\nu}=\frac{\partial v}{\partial\nu}=0,
	\qquad & x\in\pO, \ t>0, \\[1mm]
	& & \hspace*{-13mm}
	u(x,0)=u_0(x),
	\qquad & x\in\Omega, 
 	\end{array} \right.
\ee
in a bounded domain $\Omega\subset\R^n$, $n\ge 1$, where
\be{f1}
	f\in C^2([0,\infty)),
\ee
and where
\be{init}
	u_0 \in C^0(\bom)
	\quad \mbox{is nonnegative}	
	\qquad \mbox{with} \qquad
	\mint_\Omega u_0 dx=\mu>0.
\ee
Indeed, we shall see in Section \ref{sect2} that if $\Omega=B_R(0)\subset\R^n$ with some $R>0$, and if $u_0$ 
is radially symmetric, then in the resulting radial framework
the evolution can equivalently be described by considering the accumulated densities 
$w=w(s,t):=\frac{1}{n|B_1(0)|} \int_{B_r(0)} u(x,t) dx$, $s=r^n\in [0,R^n]$, $t\ge 0$, which namely solve
a Dirichlet problem for
\be{w00}
	w_t = n^2 s^{2-\frac{2}{n}} w_{ss} 
	+ n \cdot \Big(w-\frac{\mu}{n} s\Big) \cdot w_s \cdot f\Big(s^{\frac{2}{n}-2} \cdot (w-\frac{\mu}{n}s)^2 \Big)
\ee
(cf.~(\ref{0w})).
Here, to adequately quantify the destabilizing potential of the second summand on the right-hand side we shall, in a first
and yet quite basic step, rely on a comparison argument to make sure that an appropriate monotonicity assumption 
on $u_0$ entails nonnegativity of the expression $w-\frac{\mu}{n}s$ that appears in two crucial places. 
This will then enable us to suitably estimate the Burgers-type and shock-supporting action of the nonlinearity in 
(\ref{w00}) from below;
for nonlinearities which suitably generalize that in (\ref{pr}) with 
$\alpha<\frac{n-2}{2(n-1)}$ when $n\ge 3$, this will be accomplished in Section \ref{sect3} by analyzing the evolution 
of the moment-like functional 
\be{phi0}
	\int_0^{s_0} s^{-\gamma} (s_0-s) \Big(w(s,t)-\frac{\mu}{n}s\Big)  ds
\ee
along trajectories, and by thereby making sure that for smooth initial data sufficiently concentrated near the origin,
this quantity
satisfies a superlinearly forced 
ordinary differential inequality (ODI)
if the free parameter $\gamma$ herein is appropriately adjusted (Lemma \ref{lem16}).
In consequence, this will establish the following result which can be viewed as the main outcome of this study.
\begin{theo}\label{theo18}
  Let $\Omega=B_R(0) \subset \R^n$ with $n\ge 3$ and $R>0$, and let $f$ satisfy (\ref{f1}) as well as
  \be{f2}
	f(\xi) \ge \kf \cdot (1+\xi)^{-\alpha}
	\qquad \mbox{for all } \xi\ge 0
  \ee
  with some $\kf>0$ and $\alpha> 0$ fulfilling
  \be{18.1}
	\alpha<\frac{n-2}{2(n-1)}.
  \ee
  Then for any choice of $\mu>0$ one can find $R_0=R_0(\mu)\in (0,R)$ with the property that whenever 
  $u_0$ satisfies (\ref{init}) and is such that
  \be{i1}
	u_0 \mbox{ is radially symmetric}
	\qquad \mbox{with} \qquad
	\mint_{B_r(0)} u_0 dx \ge \mint_\Omega u_0 dx
	\quad \mbox{for all } r\in (0,R)
  \ee
  as well as
  \be{18.2}
	\mint_{B_{R_0}(0)} u_0 dx \ge \frac{\mu}{2} \Big(\frac{R}{R_0}\Big)^n,
  \ee
  the corresponding solution $(u,v)$ of (\ref{0})
  blows	up in finite time; that is, for the uniquely determined local classical solution, maximally extended 
  up to some time $\tm \in (0,\infty]$ according to 
  Proposition \ref{prop_loc} below,  we then have $\tm<\infty$ and
  \be{18.3}
	\limsup_{t\nearrow\tm} \|u(\cdot,t)\|_{L^\infty(\Omega)} =\infty.
  \ee
\end{theo}
{\bf Remark.} \quad
In order to construct simple examples of initial data which the above statement asserts to enforce blow-up, we only need to observe
that (\ref{i1}) is satisfied whenever $u_0$ is radial and nonincreasing with respect to $|x|$, and that (\ref{18.2})
is trivially implied if $\supp u_0 \subset \overline{B}_{R_0}(0)$, for instance.
Beyond this, however, by means of an almost verbatim copy of the reasoning detailed in \cite[Proposition 1.2]{bellomo_win_TRAN}
the set of initial data fulfilling (\ref{i1}) and (\ref{18.2}) can in fact be seen 
to contain an open subset of radial functions with respect to the topology in $L^\infty(\Omega)$,
and to furthermore even be dense in the set of all radial functions fulfilling (\ref{init})
in the framework of the topology in $L^p(\Omega)$ for each $p\in (0,1)$. 

\vspace*{1mm}
In order to secondly make sure that the above results cannot be substantially improved, 
let us next consider the case when besides (\ref{f1}), $f$ satisfies
\be{f3}
	f(\xi) \le \Kf \cdot (1+\xi)^{-\alpha} \qquad \mbox{for all } \xi\ge 0
\ee
with some $\Kf>0$ and $\alpha\in\R$. 
Then the following statement, valid even in general not necessarily radial frameworks, provides an essentially
exhausting complement to Theorem \ref{theo18}.
\begin{prop}\label{prop19}
  Let $n\ge 1$ and $\Omega\subset\R^n$ be a bounded domain with smooth boundary, and suppose that $f$ satisfies
  (\ref{f1}) and (\ref{f3}) with some $\Kf>0$ and 
  \be{19.01}
	\left\{ \begin{array}{ll}
	\alpha \in \R \qquad & \mbox{if } n=1, \\[2mm]
	\alpha > \frac{n-2}{2(n-1)} \qquad & \mbox{if } n\ge 2.
	\end{array} \right.
  \ee
  Then for any choice of $u_0$ coplying with (\ref{init}), the problem (\ref{0}) possesses a unique global classical
  solution $(u,v) \in (C^0(\bom\times [0,\infty)) \cap C^{2,1}(\bom\times (0,\infty))) \cap C^{2,0}(\bom\times (0,\infty))$
  which is bounded in the sense that
  \be{19.1}
	\|u(\cdot,t)\|_{L^\infty(\Omega)} \le C
	\qquad \mbox{for all } t>0
  \ee
  with some $C>0$.
\end{prop}
A result similar to that of Proposition \ref{prop19} has been the objective of a previous study (\cite{NT})
in which boundedness statements have been derived for (\ref{0}) in the particular case when $f(\xi)=\xi^{-\alpha}$,
$\xi> 0$, within the slightly restricted range determined by
\bas
	\left\{ \begin{array}{ll}
	\alpha \in (-\infty,\frac{1}{2}) \qquad & \mbox{if } n=1, \\[2mm]
	\alpha \in \Big( \frac{n-2}{2(n-1)}, \frac{1}{2}\Big) \qquad & \mbox{if } n\ge 2.
	\end{array} \right.
\eas
Our approach toward Proposition \ref{prop19}, essentially reducing to
a one-page argument presented in Section \ref{sect4}, apparently provides a somewhat shorter reasoning which, as we remark here
without pursuing details in this regard, can readily be extended so as to cover this result as well.\abs
We have to leave open here the question how far the number $\alpha=\frac{n-2}{2(n-1)}$, thus playing the role of a critical exponent
in the style of the definition in (\ref{ac}), belongs to either the blow-up range or the opposite regime. 
In view of precedents concerned with, e.g., the case $f\equiv 1$ when $n=2$, it may well be conceivable that also 
for arbitrary $n\ge 3$, choosing $f$ as in (\ref{pr}) with $\alpha=\frac{n-2}{2(n-1)}$  may enforce critical mass phenomena
with respect to finite-time blow-up, and possible even go along with some global unbounded solutions;
a refined analysis in this direction, however, would go beyond the scope of this study.
\mysection{Local existence and transformation to a scalar problem}\label{sect2}
In a straightforward manner adopting arguments well-established in the context of parabolic-elliptic 
Keller-Segel type systems (see e.g.~\cite{nagai1995}, \cite{djie_win} and \cite{cieslak_win} for suitable precedents), 
one can readily establish the following statement on local existence and extensibility
of solutions to (\ref{0}):
\begin{prop}\label{prop_loc}
  Let $n\ge 1$ and $\Omega\subset\R^n$ be a bounded domain with smooth boundary, 
  and assume that $f$ and $u_0$ satisfy (\ref{f1}) and (\ref{init}).
  Then there exist $\tm\in (0,\infty]$ and a uniquely determined pair $(u,v)$ of functions
  \bas
	\left\{ \begin{array}{l}
	u \in C^0(\bom\times [0,\tm)) \cap C^{2,1}(\bom\times (0,\tm)), \\[1mm]
	v\in \bigcap_{q>n} L^\infty_{loc}([0,\tm);W^{1,q}(\Omega)) \cap C^{2,0}(\bom\times (0,\tm)),
	\end{array} \right.
  \eas
  with $u\ge 0$ and $v\ge 0$ in $\Omega\times (0,\tm)$, such that 
  $(u,v)$ solves (\ref{0}) classically in $\Omega\times (0,\tm)$, that
  \be{mass}
	\io u(\cdot,t) = \io u_0
	\qquad \mbox{for all } t\in (0,\tm),
  \ee
  and that
  \be{ext}
	\mbox{if $\tm<\infty$, \quad then \quad }
	\limsup_{t\nearrow\tm} \|u(\cdot,t)\|_{L^\infty(\Omega)} = \infty.
  \ee
  Moreover, if $\Omega=B_R(0)$ with some $R>0$ and $u_0$ is radially symmetric with respect to $x=0$, then also
  $u(\cdot,t)$ and $v(\cdot,t)$ are radially symmetric for each $t\in (0,\tm)$.
\end{prop}
\mysection{Blow-up of radial solutions when $\alpha<\frac{n-2}{2(n-1)}$. Proof of Theorem \ref{theo18}}\label{sect3}
\subsection{A basic differential inequality for a moment-like functional $\phi$}
Throughout this section assuming that $\Omega=B_R(0)\subset\R^n$ is a ball with some $n\ge 2$ and $R>0$, 
for arbitrary $u_0=u_0(r)$ fulfilling (\ref{init}) 
we let $\tm\in (0,\infty]$ and the corresponding
radial local solution $(u,v)=(u(r,t),v(r,t))$ of (\ref{0}) be as provided by Proposition \ref{prop_loc},
and in the spirit of \cite{JL} we introduce
\be{w}
	w(s,t):=\int_0^{s^\frac{1}{n}} \rho^{n-1} u(\rho,t) d\rho, 
	\qquad s\in [0,R^n], \ t\in [0,\tm).
\ee
Then 
\be{ws}
	w_s(s,t)=\frac{1}{n} u(s^\frac{1}{n},t) \ge 0
	\qquad \mbox{for all $s\in (0,R^n)$ and } t\in (0,\tm),
\ee
and
\be{0w}
	\left\{ \begin{array}{lcll}
	w_t &=& n^2 s^{2-\frac{2}{n}} w_{ss} 
	+ n \cdot \Big(w-\frac{\mu}{n} s\Big) \cdot w_s \cdot f\Big(s^{\frac{2}{n}-2} \cdot (w-\frac{\mu}{n}s)^2 \Big),
	\qquad & s\in (0,R^n), \ t\in (0,\tm), \\[2mm]
	& & \hspace*{-14mm}
	w(0,t)=0, \quad w(R^n,t)=\frac{\mu R^n}{n},
	\qquad & t\in (0,\tm), \\[1mm]
	& & \hspace*{-14mm}
	w(s,0)=w_0(s):=\int_0^{s^\frac{1}{n}} \rho^{n-1} u_0(\rho) d\rho,
	\qquad & s\in (0,R^n).
	\end{array} \right.
\ee
The role of our extra assumption (\ref{i1}) in Theorem \ref{theo18} will then become clear through the 
following additional information.
\begin{lem}\label{lem1}
  Assume (\ref{f1}), (\ref{init}) and (\ref{i1}). Then for $w$ as in (\ref{w}) we have
  \bas
	w(s,t) \ge \frac{\mu}{n} \cdot s
	\qquad \mbox{for all $s\in (0,R^n)$ and } t\in (0,\tm).
  \eas
\end{lem}
\proof
  Since for $\uw(s,t):=\frac{\mu}{n} \cdot s$, $s\in [0,R^n], \ t\ge 0$, we have
  \bas
	\uw_t - n^2 s^{2-\frac{2}{n}} \uw_{ss} - n \cdot 
	\Big(\uw-\frac{\mu}{n} s\Big) \cdot \uw_s \cdot f\Big(s^{\frac{2}{n}-2} \cdot (\uw-\frac{\mu}{n}s)^2 \Big)
	=0
	\qquad \mbox{in } (0,R^n) \times (0,\infty)
  \eas
  with $\uw(0,t)=0$ and $\uw(R^n,t)=\frac{\mu R^n}{n}$ for all $t>0$, observing that our assumption (\ref{i1}) precisely asserts that
  \bas
	w_0(s)
	= \frac{s}{n} \mint_{B_{s^{1/n}}(0)} u_0 dx
	\ge \frac{s}{n} \cdot \mu = \uw(s,0)
	\qquad \mbox{for all } s\in (0,R^n),
  \eas
  this follows by applying 
  a comparison principle 
  (cf., e.g., \cite[Lemma 5.1]{bellomo_win_TRAN} for a version covering the present degenerate setting)
  to (\ref{0w}).
\qed
Based on the latter, namely, we can make use of a presupposed additional lower estimate of the form (\ref{f2})
for $f$ in establishing the following starting point of our subsequent blow-up analysis:
\begin{lem}\label{lem3}
  Suppose that (\ref{f1}) and (\ref{f2}) hold with some $\kf>0$ and $\alpha>0$, and let (\ref{init}) and (\ref{i1})
  be valid. 
  Then with $w$ and $w_0$ taken from (\ref{w}) and (\ref{0w}),
  \be{z}
	z(s,t):=w(s,t)-\frac{\mu}{n} \cdot s,
	\qquad s\in [0,R^n], \ t\in [0,\tm),
  \ee
  is nonnegative and satisfies
  \be{0z}
	z_t \ge  n^2 s^{2-\frac{2}{n}} z_{ss} + n\kf z \cdot \Big(1+s^{\frac{2}{n}-2} z^2 \Big)^{-\alpha} \cdot w_s
	\qquad \mbox{for all $s\in (0,R^n)$ and } t\in (0,\tm).
  \ee
\end{lem}
\proof
  The nonnegativity of $z$ has precisely been asserted by Lemma \ref{lem1}.
  Since moreover $w_s\ge 0$ by (\ref{ws}), in (\ref{0w}) we may use (\ref{f2}) to estimate
  \bas
	n \cdot \Big(w-\frac{\mu}{n}\cdot s\Big) \cdot w_s \cdot f\Big( 
	s^{\frac{2}{n}-2}
	(w-\frac{\mu}{n}s)^2\Big)
	&=& nzw_s \cdot f(s^{\frac{2}{n}-2} z^2) \\
	&\ge& nzw_s \cdot \kf (1+ s^{\frac{2}{n}-2} z^2)^{-\alpha}
  \eas
  for $s\in (0,R^n)$ and $t\in (0,\tm)$, and thus to obtain (\ref{0z}) from (\ref{0w}).
\qed
We can thereby state a basic evolution property of a moment-type functional which, unlike those considered in related 
precedents (\cite{biler}, \cite{win_NON_ct_signal_critexp}), explicitly involves the shifted variable $z$ instead of
the accumulated density $w$ itself:
\begin{lem}\label{lem4}
  Assume (\ref{init}), (\ref{i1}), (\ref{f1}) and (\ref{f2}) with some $\kf>0$ and $\alpha>0$, and for 
  $\gamma\in (0,1)$ and $s_0\in (0,R^n)$, let
  \be{phi}
	\phi(t):=\int_0^{s_0} s^{-\gamma} (s_0-s) z(s,t) ds,
	\qquad t\in [0,\tm),
  \ee
  where $z$ is as given by (\ref{z}).
  Then $\phi\in C^0([0,\tm)) \cap C^1((0,\tm))$ with
  \bea{4.1}
	\phi'(t)
	&\ge& - n^2 \Big(2-\frac{2}{n}-\gamma) \Big(\gamma-1+\frac{2}{n}\Big) 
		\int_0^{s_0} s^{-\frac{2}{n}-\gamma} (s_0-s) z(s,t) ds \nn\\
	& & - 2n^2 \Big(2-\frac{2}{n}-\gamma\Big) \int_0^{s_0} s^{1-\frac{2}{n}-\gamma} z(s,t) ds \nn\\[2mm]
	& & + \psi(t)
	\qquad \mbox{for all } t\in (0,\tm),
  \eea
  where
  \be{psi}
	\psi(t):=n\kf 
	\int_0^{s_0} s^{-\gamma}(s_0-s) z(s,t) \cdot \Big(1+s^{\frac{2}{n}-2} z^2(s,t)\Big)^{-\alpha} \cdot w_s(s,t) ds,
	\qquad t\in (0,\tm),
  \ee
  with $w$ taken from (\ref{w}).
\end{lem}
\proof
  Since $u \in C^0(\bom\times [0,\tm))$ and $u_t\in C^0(\bom\times (0,\tm))$, a standard argument based on the dominated convergence 
  theorem ensures that indeed $\phi$ belongs to $C^0([0,\tm))$ and to $C^1((0,\tm))$, and that for all $t\in (0,\tm)$,
  \bea{4.4}
	\phi'(t)
	&=& \int_0^{s_0} s^{-\gamma} (s_0-s) z(s,t) ds \nn\\
	&\ge& n^2 \int_0^{s_0} s^{2-\frac{2}{n}-\gamma}(s_0-s) z_{ss} ds 
	+ n\kf \int_0^{s_0} s^{-\gamma}(s_0-s) z \cdot \Big(1+s^{\frac{2}{n}-2} z^2 \Big)^{-\alpha} \cdot w_s ds
  \eea
  according to (\ref{0z}). Here, two integrations by parts show that for all $t\in (0,\tm)$ we have
  \bas
	& & \hspace*{-20mm}
	n^2 \int_0^{s_0} s^{2-\frac{2}{n}-\gamma}(s_0-s) z_{ss} ds \\
	&=& - n^2 \int_0^{s_0} \Big\{ \Big(2-\frac{2}{n}-\gamma\Big) s^{1-\frac{2}{n}-\gamma} (s_0-s) - s^{2-\frac{2}{n}-\gamma}
		\Big\} z_s ds \\
	&
	\ge
	& n^2 \int_0^{s_0} \Big\{ \Big(2-\frac{2}{n}-\gamma\Big) \Big(1-\frac{2}{n}-\gamma\Big) s^{-\frac{2}{n}-\gamma} (s_0-s)
	- 2\Big(2-\frac{2}{n}-\gamma\Big) s^{1-\frac{2}{n}-\gamma} \Big\} z ds,
  \eas
  because due to the fact that $\gamma<1 \le 2-\frac{2}{n}$ and $z(\cdot,t)\in C^1([0,R^n])$ with $z(0,t)=0$
  for all $t\in (0,\tm)$, and hence $s^{1-\frac{2}{n}-\gamma}z(s,t) \to 0$ as $s\nearrow 0$ for all $t\in (0,\tm)$, we have
  \bas
	n^2 s^{2-\frac{2}{n}-\gamma} (s_0-s) z_s \bigg|_{s=0}^{s=s_0}  
	=0
	\qquad \mbox{for all } t\in (0,\tm)
  \eas
  and
  \bas
	- n^2 \cdot \Big\{ \Big( 2-\frac{2}{n}-\gamma\Big) s^{1-\frac{2}{n}-\gamma} (s_0-s) - s^{2-\frac{2}{n}-\gamma} \Big\} z
		\bigg|_{s=0}^{s=s_0}
	&=& n^2 s_0^{2-\frac{2}{n}-\gamma} z(s_0,t) \\
	&\ge& 0
	\qquad \mbox{for all } t\in (0,\tm).
  \eas
  Therefore, (\ref{4.1}) is a consequence of (\ref{4.4}).
\qed
\subsection{A first lower estimate for the cross-diffusive contribution $\psi$}
Let us next approach the core of our analysis by quantifying the cross-diffusive contribution to (\ref{4.1}) through 
a first estimate for the rightmost summand therein from below.
This will be achieved on the basis of the following elementary observation.
\begin{lem}\label{lem2}
  Let $\alpha\in\R$ and $\beta\in (0,1]$. Then
  \be{2.1}
	(1+\xi)^{-\alpha} \ge 1- \frac{\alpha_+}{\beta} \xi^\beta
	\qquad \mbox{for all } \xi\ge 0,
  \ee
  where $\sigma_+:=\max\{\sigma,0\}$ for $\sigma\in\R$.
\end{lem}
\proof
  If $\alpha\le 0$, this is obvious. In the case when $\alpha$ is positive, (\ref{2.1})
  can be verified by observing that $\varphi(\xi):=(1+\xi)^{-\alpha} - 1 + \frac{\alpha}{\beta} \xi^\beta$,
  $\xi\ge 0$, satisfies $\varphi(0)=0$ as well as $\varphi'(\xi)=-\alpha(1+\xi)^{-\alpha-1} + \alpha \xi^{\beta-1}$
  and hence $\varphi'(\xi)\ge 0$ for all $\xi>0$, because if $\xi\in (0,1)$ then $(1+\xi)^{-\alpha-1} \le 1 \le \xi^{\beta-1}$,
  while if $\xi\ge 1$ then $(1+\xi)^{-\alpha-1} \le \xi^{-1} \le \xi^{\beta-1}$.
\qed
A straightforward application thereof shows that the function $\psi$ in (\ref{psi}), up to perturbation terms to be
estimated later on after suitably fixing the artificial parameter $\beta$, 
essentially dominates an integral containing both $z_s$ and a certain power of $z$ as factors in the integrand.
\begin{lem}\label{lem5}
  Suppose that (\ref{init}), (\ref{i1}), (\ref{f1}) and (\ref{f2}) hold with some $kf>0$ and $\alpha\in\R$,
  and let $\gamma\in (0,1)$ and $s_0\in (0,R^n)$.
  Then for any choice of $\beta\in (0,1]$, the function $\psi$ from (\ref{psi}) satisfies
  \bea{5.1}
	\psi(t)
	&\ge& n\kf \int_0^{s_0} s^{(2-\frac{2}{n})\alpha-\gamma} (s_0-s) z^{1-2\alpha}(s,t) z_s(s,t) ds \nn\\
	& & - \frac{n\kf \alpha_+}{\beta} 
	\int_0^{s_0} s^{(2-\frac{2}{n})(\alpha+\beta)-\gamma} (s_0-s) z^{1-2(\alpha+\beta)}(s,t) z_s(s,t) ds \nn\\
	& & - \frac{\mu\kf \alpha_+}{\beta} 
	\int_0^{s_0} s^{(2-\frac{2}{n})(\alpha+\beta)-\gamma} (s_0-s) z^{1-2(\alpha+\beta)}(s,t) ds
	\qquad \mbox{for all } t\in (0,\tm).
  \eea
\end{lem}
\proof
  By means of Lemma \ref{lem2}, again thanks to the nonnegativity of $z$ and $w_s$ we can estimate
  \bas
	z \cdot \Big(1+ s^{\frac{2}{n}-2} z^2 \Big)^{-\alpha} \cdot w_s
	&=& 
	s^{(2-\frac{2}{n})\alpha}
	z^{1-2\alpha} \cdot \Big(1+s^{2-\frac{2}{n}} z^{-2} \Big)^{-\alpha} \cdot w_s \\
	&\ge& 
	s^{(2-\frac{2}{n})\alpha}
	z^{1-2\alpha} \cdot \Big\{ 1 - \frac{\alpha_+}{\beta} (s^{2-\frac{2}{n}} z^{-2})^\beta \Big\} \cdot w_s \\
	&=& 
	s^{(2-\frac{2}{n})\alpha}
	z^{1-2\alpha} w_s - \frac{\alpha_+}{\beta} 
	s^{(2-\frac{2}{n})(\alpha+\beta)}
 	z^{1-2(\alpha+\beta)} w_s
	\qquad \mbox{in } (0,R^n) \times (0,\tm).	
  \eas
  As $w_s=z_s + \frac{\mu}{n} \ge z_s$, this entails that
  \bas
	z \cdot \Big(1+ s^{\frac{2}{n}-2} z^2 \Big)^{-\alpha} \cdot w_s
	&\ge& 
	s^{(2-\frac{2}{n})\alpha}
	z^{1-2\alpha} z_s \\
	& & - \frac{\alpha_+}{\beta} 
	s^{(2-\frac{2}{n})(\alpha+\beta)}
	z^{1-2(\alpha+\beta)} z_s
	- \frac{\mu \alpha_+}{n\beta} 
	s^{(2-\frac{2}{n})(\alpha+\beta)}
	z^{1-2(\alpha+\beta)}
  \eas
  in $(0,R^n) \times (0,\tm)$, so that (\ref{5.1}) results from the definition (\ref{psi}) of $\psi$.
\qed
Here another integration by parts enables us to further control the first summand on the right of (\ref{5.1}) from below
by integral expressions no longer containing $z_s$, provided that $\alpha$ satisfies a condition weaker than that in
Theorem \ref{theo18}, and that the free parameter $\gamma$ lies above some threshold.
\begin{lem}\label{lem6}
  Let (\ref{f1}) and (\ref{f2}) be valid with some $\kf>0$ and $\alpha\in (-\infty,\frac{n}{2(n-1)})$, and let $\gamma\in (0,1)$
  be such that $\gamma>(2-\frac{2}{n})\alpha$.
  Then there exists $k>0$ such that whenever (\ref{init}) and (\ref{i1}) hold, with $z$ taken from (\ref{z}) and
  for any choice of $s_0\in (0,R^n)$ we have
  \bea{6.1}
	n\kf \int_0^{s_0} s^{(2-\frac{2}{n})\alpha-\gamma} (s_0-s) z^{1-2\alpha}(s,t) z_s(s,t) ds
	&\ge& k \int_0^{s_0} s^{(2-\frac{2}{n})\alpha-\gamma-1} (s_0-s) z^{2-2\alpha}(s,t) ds \nn\\
	& & + k \int_0^{s_0} s^{(2-\frac{2}{n})\alpha-\gamma} z^{2-2\alpha}(s,t) ds
  \eea
  for all $t\in (0,\tm)$.
\end{lem}
\proof
  Using that $\alpha<1$, we may integrate by parts to see that for all $t\in (0,\tm)$,
  \bea{6.2}
	\hspace*{-4mm}
	n\kf \int_0^{s_0} s^{(2-\frac{2}{n})\alpha-\gamma} (s_0-s) z^{1-2\alpha} z_s ds
	&=& \frac{n\kf}{2-2\alpha} \int_0^{s_0} s^{(2-\frac{2}{n})\alpha-\gamma} (s_0-s) (z^{2-2\alpha})_s ds \nn\\
	&=& - \frac{n\kf}{2-2\alpha} \int_0^{s_0} \partial_s \Big\{ s^{(2-\frac{2}{n})\alpha-\gamma}(s_0-s)\Big\}
		\cdot z^{2-2-\alpha} ds,
  \eea
  where we note that the corresponding boundary terms vanish again due to the fact that for each fixed $t\in (0,\tm)$,
  $z(\cdot,t)$ belongs to $C^1([0,\R^n])$ with $z(0,t)=0$ by Proposition \ref{prop_loc}:
  This, namely, implies that for any such $t$,
  \bas
	s^{(2-\frac{2}{n})\alpha-\gamma}(s_0-s) z^{2-2\alpha}(s,t)
	\le s_0 \|z_s(\cdot,t)\|_{L^\infty((0,R^n))}^{2-2\alpha} s^{(2-\frac{2}{n})\alpha-\gamma+2-2\alpha}
	\to 0
	\qquad \mbox{as } (0,s_0)\ni s \searrow 0,
  \eas
  because $(2-\frac{2}{n})\alpha-\gamma+2-2\alpha=2-\frac{2\alpha}{n}-\gamma \ge 2-\frac{2}{n}-\gamma>0$
  as a consequence of the inequality 
  $\gamma<1\le 2-\frac{2}{n}$.
  Furthermore computing
  \bas
	\partial_s \Big\{ s^{(2-\frac{2}{n})\alpha-\gamma}(s_0-s)\Big\}
	= -\Big[ \gamma-\Big(2-\frac{2}{n}\Big)\alpha \Big] s^{(2-\frac{2}{n})\alpha-\gamma-1}(s_0-s)
	- s^{(2-\frac{2}{n})\alpha-\gamma},
	\qquad s\in (0,s_0),
  \eas
  we readily infer (\ref{6.1}) from (\ref{6.2}) with
  $k:=
  \frac{n\kf}{2-2\alpha}
  \cdot \min \{ \gamma-(2-\frac{2}{n})\alpha \ , \ 1\}$ being positive since
  $\gamma>(2-\frac{2}{n})\alpha$.
\qed
\subsection{Controlling the ill-signed summands in (\ref{5.1}). A refined lower estimate for $\psi$}
When next estimating the second and third summands on the right of (\ref{5.1}), we may evidently concentrate on the case when $\alpha$ 
is positive. 
By imposing a suitable smallness condition on the auxiliary parameter $\beta$, we may 
first rewrite the first of the respective integrals through another integration by parts.
\begin{lem}\label{lem8}
  Let (\ref{f1}) and (\ref{f2}) be valid with some $\kf>0$ and $\alpha\in (0,\frac{n}{2(n-1)})$, and let $\gamma\in (0,1)$.
  Then for any $\beta>0$ such that $\beta<1-\alpha$ and each $s_0\in (0,R^n)$, with $z$ as in (\ref{z}) we have
  \bea{8.1}
	& & \hspace*{-28mm}
	\int_0^{s_0} s^{(2-\frac{2}{n})(\alpha+\beta)-\gamma} (s_0-s) z^{1-2(\alpha+\beta)}(s,t) z_s(s,t) ds \nn\\
	&=& \frac{\gamma-(2-\frac{2}{n})(\alpha+\beta)}{2-2(\alpha+\beta)}
	\int_0^{s_0} s^{(2-\frac{2}{n})(\alpha+\beta)-\gamma-1} (s_0-s) z^{2-2(\alpha+\beta)}(s,t) ds \nn\\
	& & + \frac{1}{2-2(\alpha+\beta)}
	\int_0^{s_0} s^{(2-\frac{2}{n})(\alpha+\beta)-\gamma} z^{2-2(\alpha+\beta)}(s,t) ds 
  \eea
  for all $t\in (0,\tm)$.
\end{lem}
\proof
  Since $\alpha+\beta<1$, once more integrating by parts we compute
  \be{8.2}
	\int_0^{s_0} s^{(2-\frac{2}{n})(\alpha+\beta)-\gamma} (s_0-s) z^{1-2(\alpha+\beta)} z_s ds
	= - \frac{1}{2-2(\alpha+\beta)} 
	\int_0^{s_0} \partial_s \Big\{ s^{(2-\frac{2}{n})(\alpha+\beta)-\gamma} (s_0-s)\Big\} \cdot z^{2-2(\alpha+\beta)} ds
  \ee
  for all $t\in (0,\tm)$, where again no additional boundary terms appear, because
  \bas
	s^{(2-\frac{2}{n})(\alpha+\beta)-\gamma} (s_0-s) z^{2-2(\alpha+\beta)} (s,t)
	\le s_0 \|z_s(\cdot,t)\|_{L^\infty((0,R^n)}^{2-2(\alpha+\beta)} s^{2-\frac{2}{n}(\alpha+\beta)-\gamma} \to 0
	\qquad \mbox{as } (0,s_0)\ni s\searrow 0
  \eas
  due to the fact that $\gamma<1<2-\frac{2}{n}<2-\frac{2}{n}(\alpha+\beta)$.
  Evaluating the right-hand side in (\ref{8.2}) directly yields (\ref{8.1}).
\qed
By means of Young's inequality, we can make sure that both integrals on the right of (\ref{8.1}) can appropriately be absorbed
by expressions of the form in (\ref{6.1}), up to addition of some error terms merely depending on the potentially small
parameter $s_0$.
\begin{lem}\label{lem9}
  Let (\ref{f1}) and (\ref{f2}) be valid with some $\kf>0$ and $\alpha\in (0,\frac{n}{2(n-1)})$, and let $\gamma\in (0,1)$.
  Then given any $\beta>0$ such that $\beta<1-\alpha$, for each $\eps>0$ one can find $C=C(\beta,\eps)>0$
  such that whenever (\ref{init}) and (\ref{i1}) hold and $s_0\in (0,R^n)$,
  with $z$ 
  taken from (\ref{z}) 
  we have
  \bea{9.1}
	& & \hspace*{-20mm}
	\int_0^{s_0} s^{(2-\frac{2}{n})(\alpha+\beta)-\gamma-1} (s_0-s) z^{2-2(\alpha+\beta)}(s,t) ds \nn\\
	&\le& \eps \int_0^{s_0} s^{(2-\frac{2}{n})\alpha-\gamma-1} (s_0-s) z^{2-2\alpha}(s,t) ds 
	+ C s_0^{3-\frac{2}{n}-\gamma}
	\qquad \mbox{for all } t\in (0,\tm)
  \eea
  as well as
  \bea{9.2}
	& & \hspace*{-20mm}
	\int_0^{s_0} s^{(2-\frac{2}{n})(\alpha+\beta)-\gamma} z^{2-2(\alpha+\beta)}(s,t) ds \nn\\
	&\le& \eps \int_0^{s_0} s^{(2-\frac{2}{n})\alpha-\gamma} z^{2-2\alpha}(s,t) ds 
	+ C s_0^{3-\frac{2}{n}-\gamma}
	\qquad \mbox{for all } t\in (0,\tm)
  \eea
\end{lem}
\proof
  According to our assumption that $\alpha+\beta<1$, we may use Young's inequality to find $c_1=c_1(\beta,\eps)>0$ such that
  \bea{9.3}
	& & \hspace*{-10mm}
	\int_0^{s_0} s^{(2-\frac{2}{n})(\alpha+\beta)-\gamma-1} (s_0-s) z^{2-2(\alpha+\beta)} ds \nn\\
	&=& \int_0^{s_0} \Big\{ s^{(2-\frac{2}{n})\alpha-\gamma-1}(s_0-s) z^{2-2\alpha} \Big\}^\frac{1-\alpha-\beta}{1-\alpha}
	\cdot \Big\{ 
	s^{(2-\frac{2}{n})(\alpha+\beta)-\gamma-1 + \frac{[\gamma+1-(2-\frac{2}{n})\alpha] \cdot (1-\alpha-\beta)}{1-\alpha}}
	\cdot (s_0-s)^\frac{\beta}{1-\alpha} \Big\} ds \nn\\
	&\le& \eps \int_0^{s_0} s^{(2-\frac{2}{n})\alpha-\gamma-1} (s_0-s) z^{2-2\alpha} ds \nn\\
	& & + c_1 \int_0^{s_0} 
	s^\frac{[(2-\frac{2}{n})(\alpha+\beta)-\gamma-1] \cdot (1-\alpha)
		+ [\gamma+1-(2-\frac{2}{n})\alpha] \cdot (1-\alpha-\beta)}{\beta}
	\cdot (s_0-s) ds
	\qquad \mbox{for all } t\in (0,\tm),
  \eea
  where simplifying in the last integrand shows that
  \bas
	\int_0^{s_0} 
	s^\frac{[(2-\frac{2}{n})(\alpha+\beta)-\gamma-1] \cdot (1-\alpha)
		+ [\gamma+1-(2-\frac{2}{n})\alpha] \cdot (1-\alpha-\beta)}{\beta}
	\cdot (s_0-s) ds
	&=& \int_0^{s_0} s^{1-\frac{2}{n}-\gamma} (s_0-s) ds \\
	&\le& s_0 \int_0^{s_0} s^{1-\frac{2}{n}-\gamma} ds
	= \frac{s_0^{3-\frac{2}{n}-\gamma}}{2-\frac{2}{n}-\gamma},
  \eas
  because $\gamma<1<2-\frac{2}{n}$.
  Therefore, (\ref{9.3}) implies (\ref{9.1}), whereas (\ref{9.2}) can similarly be derived by once more using
  Young's inequality to obtain $c_2=c_2(\beta,\eps)>0$ fulfilling
  \bas
	\int_0^{s_0} s^{(2-\frac{2}{n})(\alpha+\beta)-\gamma} z^{2-2(\alpha+\beta)} ds
	&=& \int_0^{s_0} \Big\{ s^{(2-\frac{2}{n})\alpha-\gamma} z^{2-2\alpha} \Big\}^\frac{1-\alpha-\beta}{1-\alpha}
	\cdot \Big\{ 
	s^{(2-\frac{2}{n})(\alpha+\beta)-\gamma + \frac{[\gamma-(2-\frac{2}{n})\alpha] \cdot (1-\alpha-\beta)}{1-\alpha}} 
	\Big\} ds\\
	&\le& \eps \int_0^{s_0} s^{(2-\frac{2}{n})\alpha-\gamma} z^{2-2\alpha} ds \nn\\
	& & + c_2 \int_0^{s_0} 
	s^\frac{[(2-\frac{2}{n})(\alpha+\beta)-\gamma] \cdot (1-\alpha)
		+ [\gamma-(2-\frac{2}{n})\alpha] \cdot (1-\alpha-\beta)}{\beta} ds \\
	&=& \eps \int_0^{s_0} s^{(2-\frac{2}{n})\alpha-\gamma} z^{2-2\alpha} ds 
	+ c_2 \int_0^{s_0} s^{2-\frac{2}{n}-\gamma} ds \\
	&=& \eps \int_0^{s_0} s^{(2-\frac{2}{n})\alpha-\gamma} z^{2-2\alpha} ds
	+ \frac{c_2 s_0^{3-\frac{2}{n}-\gamma}}{3-\frac{2}{n}-\gamma} 
  \eas
  for all $t\in (0,\tm)$ and any choice of $s_0\in (0,R^n)$.
\qed
The last integral on the right of (\ref{5.1}) can directly be estimated in quite a similar manner.
\begin{lem}\label{lem11}
  Let (\ref{f1}) and (\ref{f2}) be valid with some $\kf>0$ and 
  $\alpha\in (0,\frac{1}{2})$,
  and let $\gamma\in (0,1)$ and $\beta>0$ satisfy $2(\alpha+ \beta)<1$.
  Then for all $\eps>0$ there exists $C=C(\beta,\eps)>0$
  such that if (\ref{init}) and (\ref{i1}) are valid and $s_0\in (0,R^n)$,
  then the function $z$ from (\ref{z}) has the property that
  \bea{11.1}
	& & \hspace*{-20mm}
	\int_0^{s_0} s^{(2-\frac{2}{n})(\alpha+\beta)-\gamma} (s_0-s) z^{1-2(\alpha+\beta)}(s,t) ds \nn\\
	&\le& \eps \int_0^{s_0} s^{(2-\frac{2}{n})\alpha-\gamma-1} (s_0-s) z^{2-2\alpha}(s,t) ds 
	+ C s_0^\frac{3-\frac{2\alpha}{n}-\gamma+2\beta(3-\frac{2}{n}-\gamma)}{1+2\beta}
	\qquad \mbox{for all } t\in (0,\tm).
  \eea
\end{lem}
\proof
  Using that $1-2(\alpha+\beta)$ is positive, again by means of Young's inequality we can find $c_1=c_1(\beta,\eps)>0$ 
  such that
  \bea{11.2}
	& & \hspace*{-10mm}
	\int_0^{s_0} s^{(2-\frac{2}{n})(\alpha+\beta)-\gamma} (s_0-s) z^{1-2(\alpha+\beta)} ds \nn\\
	&=& \int_0^{s_0} \Big\{ s^{(2-\frac{2}{n})\alpha-\gamma-1}(s_0-s) z^{2-2\alpha} 
		\Big\}^\frac{1-2(\alpha+\beta)}{2-2\alpha}
	\cdot \Big\{ 
	s^{(2-\frac{2}{n})(\alpha+\beta)-\gamma + \frac{[\gamma+1-(2-\frac{2}{n})\alpha] \cdot [1-2(\alpha+\beta)]}
		{2-2\alpha}}
	\cdot (s_0-s)^\frac{1+2\beta}{2-2\alpha} \Big\} ds \nn\\
	&\le& \eps \int_0^{s_0} s^{(2-\frac{2}{n})\alpha-\gamma-1} (s_0-s) z^{2-2\alpha} ds \nn\\
	& & + c_1 \int_0^{s_0} 
	s^\frac{[(2-\frac{2}{n})(\alpha+\beta)-\gamma] \cdot (2-2\alpha)
		+ [\gamma+1-(2-\frac{2}{n})\alpha] \cdot [1-2(\alpha+\beta)]}{1+2\beta}
	\cdot (s_0-s) ds
	\qquad \mbox{for all } t\in (0,\tm).
  \eea
  Here,
  \bas
	& & \hspace*{-40mm}
	\int_0^{s_0} 
	s^\frac{[(2-\frac{2}{n})(\alpha+\beta)-\gamma] \cdot (2-2\alpha)
		+ [\gamma+1-(2-\frac{2}{n})\alpha] \cdot [1-2(\alpha+\beta)]}{1+2\beta}
	\cdot (s_0-s) ds \\
	&=& \int_0^{s_0} s^\frac{1-\frac{2\alpha}{n}-\gamma+2\beta(1-\frac{2}{n}-\gamma)}{1+2\beta} \cdot (s_0-s) ds \\
	&\le& s_0 \int_0^{s_0} s^\frac{1-\frac{2\alpha}{n}-\gamma+2\beta(1-\frac{2}{n}-\gamma)}{1+2\beta} ds\\
	&=& \frac{1+2\beta}{2-\frac{2\alpha}{n}-\gamma+2\beta(2-\frac{2}{n}-\gamma)}
	\cdot s_0^\frac{3-\frac{2\alpha}{n}-\gamma+2\beta(3-\frac{2}{n}-\gamma)}{1+2\beta}
  \eas
  due to the fact that since $\beta>0$, $\alpha<1$ and $\gamma<1<2-\frac{2}{n}$,
  \bas
	2-\frac{2\alpha}{n}-\gamma+2\beta \Big(2-\frac{2}{n}-\gamma\Big)
	> 2-\frac{2\alpha}{n}-\gamma>2-\frac{2}{n}-\gamma>0.
  \eas
  Therefore, (\ref{11.1}) results from (\ref{11.2}).
\qed
We now combine Lemma \ref{lem8} with Lemma \ref{lem9} and Lemma \ref{lem11} to see upon suitably fixing $\beta$ that
Lemma \ref{lem6} entails the following refined lower estimate for the function $\psi$ in (\ref{psi}).
\begin{lem}\label{lem12}
  Let (\ref{f1}) and (\ref{f2}) be valid with some $\kf>0$ and 
  $\alpha\in (-\infty,\frac{1}{2})$,
  and let $\gamma\in (0,1)$ be such that $\gamma>(2-\frac{2}{n})\alpha$.
  Then one can find $C>0$ with 
  the
  property that whenever (\ref{init}) and (\ref{i1}) hold, 
  for any choice of $s_0\in (0,R^n)$ the functions $\psi$ and $z$ from (\ref{psi}) and (\ref{z}) satisfy
  \bea{12.1}
	\psi(t)
	&\ge& \frac{k}{2} \int_0^{s_0} s^{(2-\frac{2}{n})\alpha-\gamma-1} (s_0-s) z^{2-2\alpha}(s,t) ds \nn\\
	& & + \frac{k}{2} \int_0^{s_0} s^{(2-\frac{2}{n})\alpha-\gamma} z^{2-2\alpha}(s,t) ds  \nn\\
	& & - C s_0^{3-\frac{2}{n}-\gamma}
	\qquad \mbox{for all } t\in (0,\tm),
  \eea
  where $k>0$ is as given by Lemma \ref{lem6}.
\end{lem}
\proof
  In view of Lemma \ref{lem5} and Lemma \ref{lem6}, we only need to consider the case when $\alpha$ is positive, in which
  we use that the restriction $\alpha < \frac{n}{2(n-1)}$ particularly requires that $\alpha<\frac{1}{2}$, 
  whence we can fix $\beta>0$ such that $2(\alpha+\beta)<1$.
  An application of Lemma \ref{lem9} to suitably small $\eps>0$ then yields $c_1>0$ and $c_2>0$ such that
  \bas
	& & \hspace*{-20mm}
	\frac{n\kf \alpha}{\beta} \cdot \frac{\gamma-(2-\frac{2}{n})(\alpha+\beta)}{2-2(\alpha+\beta)}
	\int_0^{s_0} s^{(2-\frac{2}{n})(\alpha+\beta)-\gamma-1}(s_0-s) z^{2-2(\alpha+\beta)} ds \\
	&\le& \frac{k}{4} \int_0^{s_0} s^{(2-\frac{2}{n})\alpha-\gamma-1}(s_0-s) z^{2-2\alpha} ds
	+ c_1 s_0^{3-\frac{2}{n}-\gamma}
	\qquad \mbox{for all } t\in (0,\tm)
  \eas
  and
  \bas
	& & \hspace*{-20mm}
	\frac{n\kf \alpha}{\beta} \cdot \frac{1}{2-2(\alpha+\beta)}
	\int_0^{s_0} s^{(2-\frac{2}{n})(\alpha+\beta)-\gamma} z^{2-2(\alpha+\beta)} ds \\
	&\le& \frac{k}{2} \int_0^{s_0} s^{(2-\frac{2}{n})\alpha-\gamma} z^{2-2\alpha} ds
	+ c_2 s_0^{3-\frac{2}{n}-\gamma}
	\qquad \mbox{for all } t\in (0,\tm),
  \eas
  while using Lemma \ref{lem11} we similarly find $c_3>0$ fulfilling
  \bas
	& & \hspace*{-20mm}
	\frac{\mu \kf \alpha}{\beta} 
	\int_0^{s_0} s^{(2-\frac{2}{n})(\alpha+\beta)-\gamma}(s_0-s) z^{1-2(\alpha+\beta)} ds \\
	&\le& \frac{k}{4} \int_0^{s_0} s^{(2-\frac{2}{n})\alpha-\gamma-1}(s_0-s) z^{2-2\alpha} ds
	+ c_3 s_0^\frac{3-\frac{2\alpha}{n}-\gamma+2\beta(3-\frac{2}{n}-\gamma)}{1+2\beta}
	\qquad \mbox{for all } t\in (0,\tm).
  \eas
  In light of Lemma \ref{lem5}, Lemma \ref{lem6} and Lemma \ref{lem8}, combining these inequalities shows that
  \bea{12.2}
	\psi(t)
	&\ge& \Big(k-\frac{k}{4}-\frac{k}{4}\Big)
	\int_0^{s_0} s^{(2-\frac{2}{n})\alpha-\gamma-1}(s_0-s) z^{2-2\alpha} ds
	+\Big(k-\frac{k}{2}\Big) 
	\int_0^{s_0} s^{(2-\frac{2}{n})\alpha-\gamma} z^{2-2\alpha} ds \nn\\
	& & - (c_1+c_2) s_0^{3-\frac{2}{n}-\gamma}
	- c_3 s_0^\frac{3-\frac{2\alpha}{n}-\gamma+2\beta(3-\frac{2}{n}-\gamma)}{1+2\beta}
	\qquad \mbox{for all } t\in (0,\tm).
  \eea
  Since
  \bas
	\frac{3-\frac{2\alpha}{n}-\gamma+2\beta(3-\frac{2}{n}-\gamma)}{1+2\beta} - \Big(3-\frac{2}{n}-\gamma\Big)
	=\frac{2(1-\alpha)}{n(1+2\beta)}
  \eas
  is positive, and since thus 
  \bas
	c_3 s_0^\frac{3-\frac{2\alpha}{n}-\gamma+2\beta(3-\frac{2}{n}-\gamma)}{1+2\beta}
	\le c_3 R^\frac{2(1-\alpha)}{1+2\beta} s_0^{3-\frac{2}{n}-\gamma}
  \eas
  for any $s_0\in (0,R^n)$, from (\ref{12.2}) we directly obtain (\ref{12.1}).
\qed
\subsection{Estimating the first two integrals on the right of (\ref{4.1})}
Let us next examine how far also the two first integrals appearing on the right-hand side of (\ref{4.1}) can 
be controlled by $\psi$. Indeed, further applications of Young's inequality show that each of these integrals can essentially
be estimated against one of the first two summands on the right of (\ref{12.1}), provided that $\gamma$ satisfies an additional
smallness condition
which can be fulfilled within a range of $\alpha$ which is yet larger than that indicated in Theorem \ref{theo18}.
\begin{lem}\label{lem13}
  Suppose that (\ref{f1}) and (\ref{f2}) hold with some $\kf>0$ and 
  $\alpha\in (-\infty,\frac{n-2}{2n-3})$,
  and assume that $\gamma\in (0,1)$ is such that 
  \be{13.1}
	(1-2\alpha)\gamma < 2-\frac{4}{n} - 4\alpha + \frac{6\alpha}{n}.
  \ee
  Then for all $\eps>0$ there exists $C=C(\eps)>0$ such that if (\ref{init}) and (\ref{i1}) are valid, then 
  for arbitrary $s_0\in (0,R^n)$ the function $z$ in (\ref{z}) satisfies
  \bea{13.2}
	\int_0^{s_0} s^{-\frac{2}{n}-\gamma} (s_0-s) z(s,t) ds 
	&\le& \eps
	\int_0^{s_0} s^{(2-\frac{2}{n})\alpha-\gamma-1} (s_0-s) z^{2-2\alpha}(s,t) ds \nn\\
	& & + Cs_0^\frac{3-\frac{4}{n}-6\alpha+\frac{6\alpha}{n}-(1-2\alpha)\gamma}{1-2\alpha}
	\qquad \mbox{for all } t\in (0,\tm)
  \eea
  and
  \bea{13.3}
	\int_0^{s_0} s^{1-\frac{2}{n}-\gamma} z(s,t) ds 
	&\le& \eps
	\int_0^{s_0} s^{(2-\frac{2}{n})\alpha
	-\gamma
	} z^{2-2\alpha}(s,t) ds \nn\\
	& & + Cs_0^\frac{3-\frac{4}{n}-6\alpha+\frac{6\alpha}{n}-(1-2\alpha)\gamma}{1-2\alpha}
	\qquad \mbox{for all } t\in (0,\tm).
  \eea
\end{lem}
\proof
  Based on the fact that $2-2\alpha>1$, we once again invoke Young's inequality and thereby find $c_1=c_1(\eps)>0$ such that
  for all $t\in (0,\tm)$,
  \bea{13.4}
	& & \hspace*{-20mm}
	\int_0^{s_0} s^{-\frac{2}{n}-\gamma} (s_0-s) z ds \nn\\
	&=& \int_0^{s_0} \Big\{ s^{(2-\frac{2}{n})\alpha-\gamma-1} (s_0-s) z^{2-2\alpha} \Big\}^\frac{1}{2-2\alpha} \cdot
	\Big\{ s^{-\frac{2}{n}-\gamma + \frac{\gamma+1-(2-\frac{2}{n})\alpha}{2-2\alpha}} 
		(s_0-s)^\frac{1-2\alpha}{2-2\alpha} \Big\} ds \nn\\
	&\le& \eps 
	\int_0^{s_0} s^{(2-\frac{2}{n})\alpha-\gamma-1} (s_0-s) z^{2-2\alpha} ds
	+ c_1 \int_0^{s_0} s^\frac{(-\frac{2}{n}-\gamma)(2-2\alpha) + \gamma+1-(2-\frac{2}{n})\alpha}{1-2\alpha} (s_0-s) ds,
  \eea
  and note that here the assumption (\ref{13.1}) warrants that the exponent in the latter integral satisfies
  \bas
	\kappa
	&:=& \frac{(-\frac{2}{n}-\gamma)(2-2\alpha) + \gamma+1-(2-\frac{2}{n})\alpha}{1-2\alpha} \\
	&=& \frac{1-\frac{4}{n}-2\alpha+\frac{6\alpha}{n}-(1-2\alpha)\gamma}{1-2\alpha} \\
	&>& \frac{1-\frac{4}{n}-2\alpha+\frac{6\alpha}{n}-\Big\{ 2-\frac{4}{n}-4\alpha+\frac{6\alpha}{n}\Big\}}{1-2\alpha} \\
	&=& -1.
  \eas
  Therefore,
  \bas
	\int_0^{s_0} s^\frac{(-\frac{2}{n}-\gamma)(2-2\alpha) + \gamma+1-(2-\frac{2}{n})\alpha}{1-2\alpha} (s_0-s) ds
	&\le& s_0 \int_0^{s_0} s^\kappa ds \\
	&=& \frac{1}{1+\kappa} \cdot s_0^{2+\kappa} \\
	&=& \frac{1}{1+\kappa} \cdot s_0^\frac{3-\frac{4}{n}-6\alpha+\frac{6\alpha}{n}-(1-2\alpha)\gamma}{1-2\alpha},
  \eas
  whence (\ref{13.4}) implies (\ref{13.2}).\\
  Likewise, Young's inequality provides $c_2=c_2(\eps)>0$ such that
  \bas
	\int_0^{s_0} s^{1-\frac{2}{n}-\gamma} z ds
	&=& \int_0^{s_0} \Big\{ s^{(2-\frac{2}{n})\alpha-\gamma} z^{2-2\alpha} \Big\}^\frac{1}{2-2\alpha} \cdot
	\Big\{ s^{1-\frac{2}{n}-\gamma + \frac{\gamma-(2-\frac{2}{n})\alpha}{2-2\alpha}}  \Big\} ds \\
	&\le& \eps 
	\int_0^{s_0} s^{(2-\frac{2}{n})\alpha-\gamma} z^{2-2\alpha} ds
	+ c_2 \int_0^{s_0} s^\frac{(1-\frac{2}{n}-\gamma)(2-2\alpha) + \gamma-(2-\frac{2}{n})\alpha}{1-2\alpha} ds \\
	&=& \eps 
	\int_0^{s_0} s^{(2-\frac{2}{n})\alpha-\gamma} z^{2-2\alpha} ds
	+ c_2 \int_0^{s_0} s^\frac{2-\frac{4}{n}-4\alpha+\frac{6\alpha}{n}-(1-2\alpha)\gamma}{1-2\alpha} ds \\
	&=& \eps 
	\int_0^{s_0} s^{(2-\frac{2}{n})\alpha-\gamma} z^{2-2\alpha} ds \\
	& & + c_2 \cdot \frac{1-2\alpha}{3-\frac{4}{n}-6\alpha+\frac{6\alpha}{n}-(1-2\alpha)\gamma} \cdot
	s_0^\frac{3-\frac{4}{n}-6\alpha+\frac{6\alpha}{n}-(1-2\alpha)\gamma}{1-2\alpha}
	\qquad \mbox{for all } t\in (0,\tm),
  \eas
  and thus establishes (\ref{13.3}).
\qed
\subsection{A superlinearly forced ODI for $\phi$. The role of the condition $\alpha<\frac{n-2}{2(n-1)}$}
Finally, if $\gamma$ satisfies one further smallness condition, then the first integral on the right of (\ref{12.1})
can be identified as a potentially explosion-enforcing contribution: 
Namely, for such $\gamma$ the following lemma reveals that the superlinear and actually convex growth of 
$0\le \widehat{z} \mapsto \widehat{z}^{2-2\alpha}$ allows for an estimation of said integral in terms of a superlinear
expression of $\phi$:
\begin{lem}\label{lem14}
  Assume (\ref{f1}) and (\ref{f2}) with constants $\kf>0$ and $\alpha\in (-\infty,\frac{1}{2})$,
  and let $\gamma\in (0,1)$ satisfy
  \be{14.1}
	(1-2\alpha)\gamma < 2- 4\alpha + \frac{2\alpha}{n}.
  \ee
  Then there exists $C>0$ such that whenever (\ref{init}) and (\ref{i1}) hold
  and $z$ and $\phi$ are as in (\ref{z}) and (\ref{phi}), we have
  \bea{14.2}
	\hspace*{-4mm}
	\int_0^{s_0} s^{(2-\frac{2}{n})\alpha-\gamma-1} (s_0-s) z^{2-2\alpha}(s,t) ds 
	\ge C s_0^{-3+6\alpha-\frac{2\alpha}{n}+(1-2\alpha)\gamma} \phi^{2-2\alpha}(t)
	\quad \mbox{for all } t\in (0,\tm).
  \eea
\end{lem}
\proof
  Again relying on the property that $2-2\alpha>1$, we may use the H\"older inequality to estimate
  \bea{14.3}
	\phi(t)
	&=& \int_0^{s_0} 
	\Big\{ s^{(2-\frac{2}{n})\alpha-\gamma-1} (s_0-s) z^{2-2\alpha} \Big\}^\frac{1}{2-2\alpha}
	\cdot
	\Big\{ s^{-\gamma+\frac{\gamma+1-(2-\frac{2}{n})\alpha}{2-2\alpha}} (s_0-s)^\frac{1-2\alpha}{2-2\alpha} \Big\} ds 
	\nn\\
	&\le& \bigg\{ \int_0^{s_0} s^{(2-\frac{2}{n})\alpha-\gamma-1} (s_0-s) z^{2-2\alpha} ds \bigg\}^\frac{1}{2-2\alpha}
	\cdot \bigg\{ 
	\int_0^{s_0} s^\frac{-(2-2\alpha)\gamma+\gamma+1-(2-\frac{2}{n})\alpha}{1-2\alpha} (s_0-s) ds 
		\bigg\}^\frac{1-2\alpha}{2-2\alpha}
  \eea
  for all $t\in (0,\tm)$,
  where using that $2-4\alpha+\frac{2\alpha}{n}-(1-2\alpha)\gamma$ is positive by (\ref{14.1}), we see that
  \bas
	\int_0^{s_0} s^\frac{-(2-2\alpha)\gamma+\gamma+1-(2-\frac{2}{n})\alpha}{1-2\alpha} (s_0-s) ds 
	&\le& s_0 \int_0^{s_0} s^\frac{-(2-2\alpha)\gamma+\gamma+1-(2-\frac{2}{n})\alpha}{1-2\alpha} ds \\
	&=& c_1 s_0^\frac{3-6\alpha+\frac{2\alpha}{n}-(1-2\alpha)\gamma}{1-2\alpha}
  \eas
  with $c_1:=\frac{1-2\alpha}{2-4\alpha+\frac{2\alpha}{n}-(1-2\alpha)\gamma}>0$.
  Rearranging (\ref{14.3}) thus shows that
  \bas
	\bigg\{ \int_0^{s_0} s^{(2-\frac{2}{n})\alpha-\gamma-1} (s_0-s) z^{2-2\alpha} ds \bigg\}^\frac{1}{2-2\alpha}
	\ge c_1^{-\frac{1-2\alpha}{2-2\alpha}} s_0^{-\frac{3-6\alpha+\frac{2\alpha}{n}-(1-2\alpha)\gamma}{2-2\alpha}} \phi(t)
	\qquad \mbox{for all } t\in (0,\tm),
  \eas
  from which (\ref{14.2}) immediately follows.
\qed
Now a crucial issue consists of clarifying how far our above, and at first glance quite technical, 
conditions on $\gamma$ can simultaneously be fulfilled. 
That this can indeed be achieved under the assumption from Theorem \ref{theo18}, essentially optimal in view of
Proposition \ref{prop19}, results from a simple consideration:
\begin{lem}\label{lem15}
  Let $n\ge 2$ and $\alpha\in\R$ be such that $\alpha<\frac{n-2}{2(n-1)}$.
  Then
  \be{15.2}
	\Big(4-\frac{4}{n}\Big) \alpha^2 - \Big(6-\frac{8}{n}\Big) \alpha + 2-\frac{4}{n}>0,
  \ee
  and there exists $\gamma \in (0,1)$ such that $\gamma>(2-\frac{2}{n})\alpha$ and that moreover
  (\ref{13.1}) as well as (\ref{14.1}) are valid.
\end{lem}
\proof
  Since
  \bas
	\Big(4-\frac{4}{n}\Big) \alpha^2 - \Big(6-\frac{8}{n}\Big) \alpha + 2-\frac{4}{n}
	=\Big(4-\frac{4}{n}\Big) \cdot \Big( \frac{n-2}{2(n-1)} - \alpha \Big) \cdot (1-\alpha),
  \eas
  it is evident that our assumption on $\alpha$ implies (\ref{15.2}).
  We can thereby ensure simultaneous availability of (\ref{13.1}) and the condition that $\gamma>(2-\frac{2}{n})\alpha$:
  Namely, since
  \bas
	\Big\{ 2-\frac{4}{n}-4\alpha+\frac{6\alpha}{n}\Big\}
	- \Big\{ \Big(2-\frac{2}{n}\Big)\alpha \cdot (1-2\alpha)\Big\}
	= \Big(4-\frac{4}{n}\Big) \alpha^2 - \Big(6-\frac{8}{n}\Big) \alpha + 2-\frac{4}{n},
  \eas
  from (\ref{15.2}) it follows that the interval 
  $J:=\Big((2-\frac{2}{n})\alpha \, , \, \frac{2-\frac{4}{n}-4\alpha+\frac{6\alpha}{n}}{1-2\alpha}\Big)$ 
  is not empty, so that noting that moreover 
  $(2-\frac{2}{n})\alpha<(2-\frac{2}{n}) \cdot \frac{n-2}{2(n-1)} = \frac{n-2}{n}<1$, we see that
  it is possible to fix $\gamma \in (0,1)$ such that $\gamma\in J$.
  Since
  \bas
	2-\frac{4}{n}-4\alpha+\frac{6\alpha}{n}
	= 2-4\alpha+\frac{2\alpha}{n} - \frac{4(1-\alpha)}{n} 
	< 2-4\alpha+\frac{2\alpha}{n}
  \eas
  due to the fact that $\alpha<1$, this furthermore entails (\ref{14.1}).
\qed
Upon choosing $\gamma$ as thus specified, we can indeed turn (\ref{4.1}) into an autonomous and superlinearly forced ODI
for $\phi$ in the following style.
\begin{lem}\label{lem16}
  Let $n\ge 2$, assume (\ref{f1}) and (\ref{f2}) with some $\kf>0$ and $\alpha\in\R$ fulfilling $\alpha<\frac{n-2}{2(n-1)}$,
  and let $\gamma \in (0,1)$ be as provided by Lemma \ref{lem15}.
  Then there exists $C>0$ such that if (\ref{init}) and (\ref{i1}) hold, for any choice of $s_0\in (0,R^n)$ the function
  $\phi$ from (\ref{phi}) satisfies
  \be{16.1}
	\phi'(t) \ge \frac{1}{C} s_0^{-3+6\alpha - \frac{2\alpha}{n}+(1-2\alpha)\gamma} \phi^{2-2\alpha}(t)
	- Cs_0^\frac{3-\frac{4}{n}-6\alpha+\frac{6\alpha}{n}-(1-2\alpha)\gamma}{1-2\alpha}
	\qquad \mbox{for all } t\in (0,\tm).
  \ee
\end{lem}
\proof
  Since $\gamma$ satisfies (\ref{13.1}), we may employ Lemma \ref{lem13} to find $c_1>0$ such that with $k>0$
  given by Lemma \ref{lem6} and $\lambda:=\frac{3-\frac{4}{n}-6\alpha+\frac{6\alpha}{n}-(1-2\alpha)\gamma}{1-2\alpha}$,
  \bas
	n^2 \Big(2-\frac{2}{n}-\gamma\Big)\Big(\gamma-1+\frac{2}{n}\Big)
	\int_0^{s_0} s^{-\frac{2}{n}-\gamma} (s_0-s) zds
	\le \frac{k}{4} \int_0^{s_0} s^{(2-\frac{2}{n})\alpha-\gamma-1} (s_0-s) z^{2-2\alpha} ds
	+ c_1 s_0^\lambda
  \eas
  and
  \bas
	n^2 \Big(2-\frac{2}{n}-\gamma\Big)
	\int_0^{s_0} s^{1-\frac{2}{n}-\gamma} zds
	\le \frac{k}{2} \int_0^{s_0} s^{(2-\frac{2}{n})\alpha-\gamma} z^{2-2\alpha} ds
	+ c_1 s_0^\lambda
  \eas
  for all $t\in (0,\tm)$.
  As furthermore $\gamma>(2-\frac{2}{n})\alpha$, we may use this in conjunction with Lemma \ref{lem4} and Lemma \ref{lem12}
  to see that with some $c_2>0$ we have
  \be{16.2}
	\phi'(t)
	 \ge \frac{k}{4} \int_0^{s_0} s^{(2-\frac{2}{n})\alpha-\gamma-1} (s_0-s) z^{2-2\alpha} ds
	- 2c_1 s_0^\lambda - c_2 s_0^{3-\frac{2}{n}-\gamma}
	\qquad \mbox{for all } t\in (0,\tm),
  \ee
  where we note that
  \bas
	\Big(3-\frac{2}{n}-\gamma\Big) -\lambda 	
	= \frac{2(1-\alpha)}{n(1-2\alpha)}
  \eas
  is positive and hence
  \be{16.3}
	c_2 s_0^{3-\frac{2}{n}-\gamma} 
	= c_2 s_0^\frac{2(1-\alpha)}{n(1-2\alpha)} \cdot s_0^\lambda
	\le c_3 s_0^\lambda
  \ee
  with $c_3:=c_2 R^\frac{2(1-\alpha)}{1-2\alpha}$.\\
  Apart from that, we may rely on the fact that $\gamma$ also satisfies (\ref{14.1}) to infer from Lemma \ref{lem14}
  the existence of $c_4>0$ such that
  \bas
	\frac{k}{4} \int_0^{s_0} s^{(2-\frac{2}{n})\alpha-\gamma-1} (s_0-s) z^{2-2\alpha} ds
	\ge c_4 s_0^{-3+6\alpha-\frac{2\alpha}{n}+(1-2\alpha)\gamma} \phi^{2-2\alpha}(t)
	\qquad \mbox{for all } t\in (0,\tm).
  \eas
  Therefore, (\ref{16.2}) together with (\ref{16.3}) establishes (\ref{16.1}) if we let
  $C:=\max\{\frac{1}{c_4} \, , \, 2c_1+c_3\}$.
\qed
\subsection{Selection of concentrated initial data. Proof of Theorem \ref{theo18}}
In order to derive a blow-up result from Lemma \ref{lem16} by means of a contradiction argument, it remains to
select the free parameter $s_0$, along with initial data which reflect appropriate mass concentration near the origin,
in such a way that (\ref{16.1}) at the initial instant, and hence throughout evolution by comparison, can be turned
into a genuine superlinear ODI without appearance of any negative summand on its right:
\begin{lem}\label{lem17}
  Let $n\ge 2$, and suppose that (\ref{f1}) and (\ref{f2}) hold with some $\kf>0$ and $\alpha\in\R$ such that
  $\alpha<\frac{n-2}{2(n-1)}$.
  Then for each $\mu>0$ one can find $s_0=s_0(\mu)\in (0,\frac{R^n}{4})$ such that whenever
  (\ref{init}) and (\ref{i1}) hold and $w_0$ from (\ref{0w}) satisfies
  \be{17.1}
	w_0 \Big(\frac{s_0}{2}\Big) \ge \frac{\mu R^n}{2n},
  \ee
  the solution of (\ref{0}) has the property that $\tm<\infty$.
\end{lem}
\proof
  Relying on our assumption $\alpha<\frac{n-2}{2(n-1)}$, we can fix $\gamma\in (0,1)$ as given by Lemma \ref{lem15},
  and then invoke Lemma \ref{lem16} to find $c_1>0$ and $c_2>0$ with the property that if (\ref{init}) and (\ref{i1}) hold,
  then for arbitrary $s_0\in (0,R^n)$, the function $\phi$ as accordingly defined through (\ref{phi}) satisfies
  \be{17.2}
	\phi'(t)
	\ge c_1 s_0^{-3+6\alpha - \frac{2\alpha}{n}+(1-2\alpha)\gamma} \phi^{2-2\alpha}(t)
	- c_2 s_0^\frac{3-\frac{4}{n}-6\alpha+\frac{6\alpha}{n}-(1-2\alpha)\gamma}{1-2\alpha}
	\qquad \mbox{for all } t\in (0,\tm).
  \ee
  We now use that $\gamma<1$ in introducing
  \bas
	c_3=c_3(\mu):=\frac{\mu R^n}{16n} \cdot \frac{(\frac{3}{4})^{1-\gamma} -(\frac{1}{2})^{1-\gamma}}{1-\gamma}>0
  \eas
  and then choose $s_0=s_0(\mu)\in (0,\frac{R^n}{4})$ small enough such that
  \be{17.3}
	\frac{1}{2} c_1 c_3^{2-2\alpha} s_0^{1+2\alpha-\frac{2\alpha}{n}
	-\gamma
	}
	\ge c_2 s_0^\frac{3-\frac{4}{n}-6\alpha+\frac{6\alpha}{n}-(1-2\alpha)\gamma}{1-2\alpha},
  \ee
  noting that the latter is possible due to (\ref{15.2}), which indeed ensures that
  \bas
	\frac{3-\frac{4}{n}-6\alpha+\frac{6\alpha}{n}-(1-2\alpha)\gamma}{1-2\alpha}
	- \Big\{ 1+2\alpha-\frac{2\alpha}{n} 
	-\gamma
	\Big\}
	= \frac{(4-\frac{4}{n})\alpha^2 -(6-\frac{8}{n})\alpha + 2-\frac{4}{n}}{1-2\alpha}
  \eas
  is positive, and that hence (\ref{17.3}) holds for all suitably small $s_0>0$.
  Keeping this value of $s_0$ fixed, we now assume to be given initial data which besides (\ref{init}) and (\ref{i1})
  fulfill (\ref{17.1}).
  Then by monotonicity of $w_0$, we have $w_0(s) \ge \frac{\mu R^n}{2n}$ for all $s\in (\frac{s_0}{2},\frac{3s_0}{4})$
  and hence, by (\ref{z}) and our definition of $c_3$,
  \bas
	\phi(0)
	&\ge& \frac{s_0}{4} \int_{\frac{s_0}{2}}^{\frac{3s_0}{4}} s^{-\gamma} z_0(s) ds \\
	&=& \frac{s_0}{4} \int_{\frac{s_0}{2}}^{\frac{3s_0}{4}} s^{-\gamma} \Big( w_0(s)-\frac{\mu}{n} \cdot s\Big) ds \\
	&\ge& \frac{\mu s_0}{4n} \int_{\frac{s_0}{2}}^{\frac{3s_0}{4}} s^{-\gamma} \Big( \frac{R^n}{2}-s\Big) ds \\
	&\ge& \frac{\mu R^n s_0}{16n} \int_{\frac{s_0}{2}}^{\frac{3s_0}{4}} s^{-\gamma} ds \\
	&=& c_3 s_0^{2-\gamma},
  \eas
  because $s_0\le \frac{R^n}{4}$.
  As a consequence of (\ref{17.3}), we thus obtain that
  \bea{17.99}
	& & \hspace*{-20mm}
	c_1 s_0^{-3+6\alpha - \frac{2\alpha}{n}+(1-2\alpha)\gamma} \phi^{2-2\alpha}(0)
	- c_2 s_0^\frac{3-\frac{4}{n}-6\alpha+\frac{6\alpha}{n}-(1-2\alpha)\gamma}{1-2\alpha} \nn\\
	&\ge& \frac{1}{2} c_1 s_0^{-3+6\alpha - \frac{2\alpha}{n}+(1-2\alpha)\gamma} \phi^{2-2\alpha}(0),
  \eea
  whence an ODE comparison argument shows that thanks to (\ref{17.2}),
  \be{17.999}
	\phi'(t) \ge \frac{1}{2} c_1 s_0^{-3+6\alpha - \frac{2\alpha}{n}+(1-2\alpha)\gamma} \phi^{2-2\alpha}(t)
	\qquad \mbox{for all } t\in (0,\tm),
  \ee
  for writing
  \bas
	\underline{\phi}(t):=\Bigg( 
	\frac{2c_2 s_0^\frac{3-\frac{4}{n}-6\alpha+\frac{6\alpha}{n}-(1-2\alpha)\gamma}{1-2\alpha}}
		{c_1 s_0^{-3+6\alpha - \frac{2\alpha}{n}+(1-2\alpha)\gamma}}
	\Bigg)^\frac{1}{2-2\alpha},
	\qquad t\ge 0,
  \eas
  we directly observe that
  \bas
	\underline{\phi}_t - c_1 s_0^{-3+6\alpha - \frac{2\alpha}{n}+(1-2\alpha)\gamma} \underline{\phi}^{2-2\alpha}(t)
	+ c_2 s_0^\frac{3-\frac{4}{n}-6\alpha+\frac{6\alpha}{n}-(1-2\alpha)\gamma}{1-2\alpha}
	&=& - c_2 s_0^\frac{3-\frac{4}{n}-6\alpha+\frac{6\alpha}{n}-(1-2\alpha)\gamma}{1-2\alpha} \\
	&\le& 0
	\qquad \mbox{for all } t>0,
  \eas
  whereas (\ref{17.99}) ensures that $\underline{\phi}(0) \le \phi(0)$.
  Since $2-2\alpha>1$, by positivity of $\phi(0)$ 
  the inequality in (\ref{17.999}), however implies 
  that $\phi$, and thus clearly also $(u,v)$, must cease to
  exist in finite time.
\qed
Our main result on the occurrence of blow-up in (\ref{0}) can now be obtained by simply transforming the above
back to the original variables:\abs
\proofc of Theorem \ref{theo18}. \quad
  With $s_0=s_0(\mu)\in (0,\frac{R^n}{4})$ as provided by Lemma \ref{lem17}, we let $R_0=R_0(\mu)
  :=(\frac{s_0}{2})^\frac{1}{n} \in (0,R)$, and assume
  $u_0$ to satisfy (\ref{init}), (\ref{i1}) as well as (\ref{18.2}).
  We then only need to observe that when rewritten in terms of the function $w_0$ defined in (\ref{0w}), 
  (\ref{18.2}) says that
  \bas
	w_0\Big(\frac{s_0}{2}\Big)
	= \frac{R_0^n}{n} \mint_{B_{R_0}(0)} u_0 dx
	\ge \frac{R_0^n}{n} \cdot \frac{\mu}{2}\Big(\frac{R}{R_0}\Big)^n
	= \frac{\mu R^n}{2n}.
  \eas
  Therefore, namely, the claim becomes a direct consequence of Lemma \ref{lem17}.
\qed
\mysection{Global boundedness for supercritical $\alpha$. Proof of Proposition \ref{prop19}}\label{sect4}
Let us finally make sure that our blow-up result indeed is essentially optimal with respect to the parameter range therein.
Our brief reasoning in this direction will rely on a basic integrability feature of the taxis gradient which immediately
results from the $L^1$ bound for $u$ implied by (\ref{mass}) due to standard elliptic regularity theory:
\begin{lem}\label{lem20}
  Under the assumptions of Proposition \ref{prop19}, for all $q\in [1,\frac{n}{n-1})$ one can find $C=C(q)>0$ such that
  \bas
	\|\nabla v(\cdot,t)\|_{L^q(\Omega)} \le C
	\qquad \mbox{for all } t\in (0,\tm).
  \eas
\end{lem}
\proof
  In view of (\ref{mass}),
  this is a direct consequence of the second equation in (\ref{0}) when combined with 
  well-known regularity theory for elliptic problems with $L^1$ inhomogeneities
  (\cite{brezis_strauss}).
\qed
We can thereby easily derive the claimed statement on global existence and boundedness by means of a suitably designed
$L^\infty$ estimation procedure based on smoothing properties of the Neumann heat semigroup:\abs
\proofc of Proposition \ref{prop19}. \quad
  Without loss of generality assuming that $\alpha<\frac{1}{2}$, we note that since $2(n-1)\alpha>n-2$ and hence
  $(1-2\alpha)n <\frac{n}{n-1}$, we can fix $q\in [1,\frac{n}{n-1})$ such that $q>(1-2\alpha)n$, whereupon it becomes
  possible to pick $r>n$ such that still $q>(1-2\alpha)r$.\\
  Then according to known smoothing properties of the Neumann heat semigroup $(e^{t\Delta})_{t\ge 0}$ in $\Omega$
  (\cite{FIWY}), there exist $c_1>0$ and $\theta>0$ such that for all $\varphi\in C^1(\bom;\R^n)$ such that
  $\varphi\cdot\nu=0$ on $\pO$,
  \bas
	\|e^{t\Delta} \nabla\cdot\varphi\|_{L^\infty(\Omega)} \le c_1 t^{-\frac{1}{2}-\frac{n}{2r}} e^{-\theta t}
		\|\varphi\|_{L^r(\Omega)}
	\qquad \mbox{for all } t>0.
  \eas
  Therefore, using 
  that $0\le e^{t\Delta} u_0 \le \|u_0\|_{L^\infty(\Omega)}$ for all $t>0$ by the maximum principle,
  we infer that for all $t\in (0,\tm)$,
  \bea{19.2}
	\hspace*{-7mm}
	\|u(\cdot,t)\|_{L^\infty(\Omega)}
	&=& \Bigg\| e^{t\Delta} u_0
	- \int_0^t e^{(t-s)\Delta} 
	\nabla \cdot \bigg\{ u(\cdot,s) f\Big( |\nabla v(\cdot,s)|^2 \Big) \nabla v(\cdot,s) 
	\bigg\} ds \Bigg\|_{L^\infty(\Omega)} \nn\\
	&\le& \|e^{t\Delta} u_0\|_{L^\infty(\Omega)}
	+ c_1 \int_0^t \bigg\| e^{(t-s)\Delta} 
	\nabla \cdot \bigg\{ u(\cdot,s) f\Big( |\nabla v(\cdot,s)|^2 \Big) \nabla v(\cdot,s) 
	\bigg\} \bigg\|_{L^\infty(\Omega)} ds \nn\\
	&\le& \|u_0\|_{L^\infty(\Omega)}
	+ c_1 \int_0^t (t-s)^{-\frac{1}{2}-\frac{n}{2r}} e^{-\theta(t-s)} 
	\bigg\| u(\cdot,s) f\Big( |\nabla v(\cdot,s)|^2 \Big) \nabla v(\cdot,s) \bigg\|_{L^r(\Omega)} ds,
  \eea
  where by (\ref{f3}) and the H\"older inequality, for all $s\in (0,\tm)$ we see that
  \bas
	\bigg\| u(\cdot,s) f\Big( |\nabla v(\cdot,s)|^2 \Big) \nabla v(\cdot,s) \bigg\|_{L^r(\Omega)}
	&\le& \Kf \Big\| u(\cdot,s) \Big(1+|\nabla v(\cdot,s)|^2 \Big)^{-\alpha} \nabla v(\cdot,s) \Big\|_{L^r(\Omega)} \\
	&\le& \Kf \Big\| u(\cdot,s) |\nabla v(\cdot,s)|^{1-2\alpha} \Big\|_{L^r(\Omega)} \\
	&\le& \Kf \|u(\cdot,s)\|_{L^\frac{qr}{q-(1-2\alpha)r}(\Omega)} \|\nabla v(\cdot,s)\|_{L^q(\Omega)}^{1-2\alpha} \\
	&\le& \Kf \|u(\cdot,s)\|_{L^\infty(\Omega)}^a \|u(\cdot,s)\|_{L^1(\Omega)}^{1-a} 
	\|\nabla v(\cdot,s)\|_{L^q(\Omega)}^{1-2\alpha}
  \eas
  with $a:=1-\frac{q-(1-2\alpha)r}{qr} \in (0,1)$.
  As $\|u(\cdot,s)\|_{L^1(\Omega)}=\io u_0$ for all $s\in (0,\tm)$,
  by means of Lemma \ref{lem20} we thus find $c_2>0$ such that writing 
  $M(T):=\sup_{t\in (0,T)} \|u(\cdot,t)\|_{L^\infty(\Omega)}$ for $T\in (0,\tm)$, we have
  \bas
	\bigg\| u(\cdot,s) f\Big( |\nabla v(\cdot,s)|^2 \Big) \nabla v(\cdot,s) \bigg\|_{L^r(\Omega)}
	\le c_2\|u(\cdot,s)\|_{L^\infty(\Omega)}^a 
	\le c_2 M^a(T)
	\qquad \mbox{for all } s\in (0,T),
  \eas
  whence (\ref{19.2}) shows that
  \bea{19.3}
	\|u(\cdot,t)\|_{L^\infty(\Omega)}
	&\le& \|u_0\|_{L^\infty(\Omega)}
	+ c_1 c_2 M^a(T) \int_0^t (t-s)^{-\frac{1}{2} - \frac{n}{2r}} e^{-\theta(t-s)} ds \nn\\
	&\le& c_3 + c_3 M^a(T)
	\qquad \mbox{for all } t\in (0,T)
  \eea
  if we let 
  $c_3:=\max\left\{ \|u_0\|_{L^\infty(\Omega)} \, , \, 
  c_1 c_2 \int_0^\infty \sigma^{-\frac{1}{2}-\frac{n}{2r}} e^{-\theta\sigma} \right\}$ and note that $c_3$ is finite due to the 
  restriction that $r>n$.
  In consequence, (\ref{19.3}) entails that
  \bas
	M(T) \le c_3 + c_3 M^a(T)
	\qquad \mbox{for all } T\in (0,\tm)
  \eas
  and thereby asserts that $\|u(\cdot,t)\|_{L^\infty(\Omega)} \le \max\{1 \, , \, (2c_3)^\frac{1}{1-a} \}$
  for all $t\in (0,\tm)$.
  Thanks to (\ref{ext}), this firstly ensures that $\tm=\infty$, and that secondly moreover (\ref{19.1}) holds.
\qed
\vspace*{5mm}
{\bf Acknowledgement.} \quad
  The author warmly thanks the anonymous reviewer for numerous helpful remarks and suggestions.
  He furthermore
  acknowledges support of the {\em Deutsche Forschungsgemeinschaft} 
  in the context of the project {\em Emergence of structures and advantages in cross-diffusion systems} 
  (Project No.~411007140, GZ: WI 3707/5-1).

\begin{thebibliography}{99}
%
\bibitem{bellomo_flim}
  \sc Bellomo, N.,  Bellouquid, A., Nieto, J., Soler, J.:
  \it Multiscale biological tissue models and flux-limited chemotaxis from binary mixtures
  of multicellular growing systems.
  \rm Math.~Mod.~Meth.~Appl.~Sci. {\bf 20},  1675-1693 (2010)
\bibitem{BBTW}
  \sc Bellomo, N., Bellouquid, A., Tao, Y., Winkler, M.:
  \it Toward a mathematical theory of Keller-Segel models of pattern formation in biological tissues. 
  \rm Math.~Models Methods Appl.~Sci. {\bf 25}, 1663-1763 (2015)
\bibitem{bellomo_win_TRAN}
  \sc Bellomo, N., Winkler, M.: \it Finite-time blow-up in a degenerate chemotaxis system with flux limitation. 
  \rm Trans.~Amer.~Math.~Soc.~Ser. B {\bf 4}, 31-67 (2017)
\bibitem{bianchi_painter_sherratt}
  \sc Bianchi, A., Painter, K.J., Sherratt, J.A.:
  \it A mathematical model for lymphangiogenesis in normal and diabetic wounds.  
  \rm J.~Theor.~Biol. {\bf 383}, 61-86 (2015)
\bibitem{bianchi_painter_sherratt2016}
  \sc Bianchi, A., Painter., K.J., Sherratt, J.A.: 
  \it Spatio-temporal models of lymphangiogenesis in wound healing. 
  \rm Bull.~Math.~Biol. {\bf 78}, 1904-1941 (2016)
\bibitem{biler}
  \sc Biler, P.:  \it Local and global solvability of some parabolic systems modelling chemotaxis. 
  \rm Adv. Math. Sci. Appl. {\bf 8}, 715-743 (1998)
\bibitem{brezis_strauss}
  \sc Br\'ezis, H., Strauss, W.A.: \it Semi-linear second-order elliptic equations in $L^1$.
  \rm J.~Math.~Soc.~Japan {\bf 25}, 565-590 (1973)
\bibitem{cieslak_laurencot_DCDS}
  \sc Cie\'slak, T., Lauren\c{c}ot, Ph.: 
  \it Looking for critical nonlinearity in the one-dimensional quasilinear Smoluchowski-Poisson system.
  \rm Discr.~Cont.~Dyn.~Syst.~A {\bf 26}, 417-430 (2010)
\bibitem{cieslak_laurencot_ANIHPC}
  \sc Cie\'slak, T., Lauren\c{c}ot, Ph.: 
  \it Finite time blow-up for a one-dimensional quasilinear parabolic-parabolic chemotaxis system. 
  \rm Ann.~Inst.~H.~Poincar\'e Anal.~Non Lin\'eaire {\bf 27} (1), 437-446 (2010)
\bibitem{cieslak_stinner_JDE2012}
  \sc Cie\'slak, T., Stinner, C.: 
  \it Finite-time blowup and global-in-time unbounded solutions to a parabolic-parabolic quasilinear Keller-Segel 
  system in higher dimensions.
  \rm J.~Differ.~Eq. {\bf 252} (10), 5832-5851 (2012)
\bibitem{cieslak_stinner_JDE2015}
  \sc Cie\'slak, T., Stinner, C.: 
  \it New critical exponents in a fully parabolic quasilinear Keller-Segel system and applications to volume filling models.
  \rm J.~Differ.~Eq. {\bf 258} (6), 2080-2113 (2015)
\bibitem{cieslak_win}
  \sc Cie\'slak, T., Winkler, M.: \it Finite-time blow-up in a quasilinear system of chemotaxis.
  \rm Nonlinearity {\bf 21}, 1057-1076 (2008)
\bibitem{deng_levine}
  \sc Deng, K., Levine, H.A.:
  \it The role of critical exponents in blow-up theorems: The sequel.
  \rm J.~Math.~Anal.~Appl. {\bf 243}, 85-126 (2000)
\bibitem{djie_win}
  \sc Djie, K., Winkler, M.: \it Boundedness and finite-time collapse in a chemotaxis system with volume-filling effect.
  \rm Nonlin.~Anal. {\bf 72}, 1044-1064 (2010)
\bibitem{fuest_logistic}
  \sc Fuest, M.: \it On finite-time blow-up in chemotaxis systems with logistic source.
  \rm Preprint
\bibitem{FIWY}
  \sc Fujie, K., Ito, A., Winkler, M., Yokota, T.:
  \it Stabilization in a chemotaxis model for tumor invasion.
  \rm Discrete Cont.~Dyn.~Syst. {\bf 36}, 151-169 (2016)
\bibitem{herrero_velazquez}
  \sc Herrero, M.A., Vel\'azquez, J.J.L.: \it A blow-up mechanism for a chemotaxis model.
  \rm Ann.~Scu.~Norm.~Super.~Pisa Cl.~Sci. {\bf 24}, 663-683 (1997)
\bibitem{hillen_painter2009}
  \sc Hillen, T., Painter, K.J.: \it A user's guide to PDE models for chemotaxis.
  \rm J.~Math.~Biol. {\bf 58} (1), 183-217 (2009)
\bibitem{horstmann_DMV}
  \sc Horstmann, D.: \it From 1970 until present: The Keller-Segel model in chemotaxis and its consequences I.
  \rm Jahresberichte DMV {\bf 105} (3), 103-165 (2003)
\bibitem{JL}
  \sc J\"ager, W., Luckhaus, S.: 
  \it On explosions of solutions to a system of partial differential equations modelling chemotaxis.
  \rm Trans.~Am.~Math.~Soc. {\bf 329}, 819-824 (1992)
\bibitem{kaplan}
  \sc Kaplan, S.: \it On the growth of solutions of quasi-linear parabolic equations.
  \rm Comm. Pure Appl. Math. {\bf 16}, 305-330 (1963)
\bibitem{KS}
  \sc Keller, E.F., Segel, L.A.: \it  Initiation of slime mold aggregation viewed as an instability.
  \rm J.~Theoret.~Biol. {\bf 26} 399-415 (1970)
\bibitem{kowalczyk_szymanska}
  \sc Kowalczyk, R., Szyma\'nska, Z.: \it On the global existence of solutions to an aggregation model.
  \rm J.~Math.~Anal.~Appl. {\bf 343}, 379-398 (2008)
\bibitem{lankeit_ITBU}
  \sc Lankeit, J.:
  \it Infinite time blow-up of many solutions to a general quasilinear parabolic-elliptic Keller-Segel system. 
  \rm Discrete Contin. Dyn. Syst., Ser. S, {\bf 13}, 233-255 (2020)
\bibitem{levine}
  \sc Levine, H.A.:
  \it The role of critical exponents in blowup theorems.
  \rm SIAM Review {\bf 32}, 262-288 (1990)
\bibitem{liu_tao}
  \sc Liu, D., Tao, Y.: \it Boundedness in a chemotaxis system with nonlinear signal production. 
  \rm Appl.~Math.~J.~Chinese Univ.~Ser.~B {\bf 31}, 379-388 (2016)
\bibitem{nagai1995}
  \sc Nagai, T.: \it Blow-up of radially symmetric solutions to a chemotaxis system.
  \rm Adv.~Math.~Sci.~Appl. {\bf 5}, 581-601 (1995)
\bibitem{nagai2001}
  \sc Nagai, T.: \it 
  Blowup of nonradial solutions to parabolic-elliptic systems modeling
  chemotaxis in two-dimensional domains.
  \rm J.~Inequal.~Appl. {\bf 6}, 37-55 (2001)
\bibitem{NSY}
  \sc Nagai, T., Senba, T., Yoshida, K.: 
  \it Application of the Trudinger-Moser inequality to a parabolic system of chemotaxis.
  \rm Funkc.~Ekvacioj, Ser.~Int. {\bf 40}, 411-433 (1997)
\bibitem{NT}
  \sc Negreanu, M., Tello, J.I.: 
  \it On a parabolic-elliptic system with gradient dependent chemotactic coefficient.
  \rm J.~Differential Eq. {\bf 265}, 733-751 (2018)
\bibitem{perthame}
  \sc Perthame, B., Yasuda, S.:
  \it Stiff-response-induced instability for chemotactic bacteria and flux-limited Keller-Segel equation.
  \rm Nonlinearity {\bf 31}, 4065 (2018)
\bibitem{senba_suzuki_AAA}
  \sc Senba, T., Suzuki, T.: \it A quasi-linear system of chemotaxis.
  \rm Abstr.~Appl.~Anal. 2006, 1-21 (2006)
\bibitem{suzuki_book}
  \sc Suzuki, T.: \it Free Energy and Self-Interacting Particles. 
  \rm Birkh\"auser, Boston (2005)
\bibitem{taowin_subcrit}
  \sc Tao, Y., Winkler, M.: \it
  Boundedness in a quasilinear parabolic-parabolic Keller-Segel system with subcritical sensitivity,
  \rm J.~Differential Eq. {\bf 252}, 692-715 (2012)
\bibitem{taowin_JEMS}
  \sc Tao, Y., Winkler, M.:
  \it Critical mass for infinite-time aggregation in a chemotaxis model with indirect signal production.
  \rm J.~European Math.~Soc. {\bf 19}, 3641-3678 (2017)
\bibitem{tello_win}
  \sc Tello, J.I., Winkler, M.: 
  \it A chemotaxis system with logistic source. 
  \rm Comm.~Part.~Differential Eq. {\bf 32} (6), 849-877 (2007)
\bibitem{win_collapse}
  \sc Winkler, M.: \it Does a `volume-filling effect' always prevent chemotactic collapse?
  \rm Math.~Meth.~Appl.~Sci. {\bf 33}, 12-24 (2010)
\bibitem{win_NON_ct_signal_critexp}
  \sc Winkler, M.: 
  \it A critical blow-up exponent in a chemotaxis system with nonlinear signal production.
  \rm Nonlinearity {\bf 31}, 2031-2056 (2018)
\bibitem{win_ITBU}
  \sc Winkler, M.: 
  \it Global classical solvability and generic infinite-time blow-up in quasilinear Keller-Segel systems 
  with bounded sensitivities.
  \rm J.~Differential Eq. {\bf 266}, 8034-8066 (2019)
%
\end{thebibliography}
\end{document}